\documentclass[12 pt]{article}
\counterwithin*{equation}{section}
\usepackage[all]{xy}
\usepackage{amsmath,amssymb,amsfonts,latexsym,float,graphics,epsfig}

\usepackage{setspace}
\usepackage[top=2cm, bottom=2cm, left=2cm, right=2cm]{geometry}
\usepackage[colorlinks]{hyperref}
\usepackage{tikz}
\usetikzlibrary{cd}
\newcommand{\contr}[1]{\iota_{#1}}
\newcommand*{\comp}{\circ}
\DeclareMathOperator{\im}{im}

\newcommand*{\RR}{\mathbb{R}}
\newcommand{\dd}{\mathrm{d}}
\newtheorem{theorem}{Theorem}[section]
\newtheorem{proposition}[theorem]{Proposition}

\newtheorem{corollary}[theorem]{Corollary}

\newtheorem{remark}{Remark}

\begin{document}
\title{Contact Dynamics versus Legendrian and Lagrangian Submanifolds} \maketitle

\begin{center}
O\u{g}ul Esen\footnote{E-mail: 
\href{mailto:oesen@gtu.edu.tr}{oesen@gtu.edu.tr}}\\
Department of Mathematics, \\ Gebze Technical University, 41400 Gebze,
Kocaeli, Turkey.

\bigskip

\author{Manuel Lainz Valcázar\footnote
{E-mail: 
\href{mailto:manuel.lainz@icmat.es}{manuel.lainz@icmat.es}}
\\ 
Instituto de Ciencias Matematicas, Campus Cantoblanco \\ 
Consejo Superior de Investigaciones Cient\'ificas
 \\
C/ Nicol\'as Cabrera, 13--15, 28049, Madrid, Spain
}

\bigskip

Manuel de Le\'on\footnote{E-mail: \href{mailto:mdeleon@icmat.es}{mdeleon@icmat.es}}
\\ Instituto de Ciencias Matem\'aticas, Campus Cantoblanco \\
 Consejo Superior de Investigaciones Cient\'ificas
 \\
C/ Nicol\'as Cabrera, 13--15, 28049, Madrid, Spain
\\
and
\\
Real Academia Espa{\~n}ola de las Ciencias.
\\
C/ Valverde, 22, 28004 Madrid, Spain.

\bigskip

Juan Carlos Marrero\footnote
{E-mail: 
\href{mailto:jcmarrer@ull.edu.es}{jcmarrer@ull.edu.es}}
\\
ULL-CSIC Geometria Diferencial y Mec\'{a}nica Geom\'etrica,\\
Departamento de Matematicas, Estadistica e I O, 
\\
Secci\'on de Matem\'aticas, 
\\
Facultad de Ciencias,
Universidad de la Laguna, La Laguna, Tenerife, Canary Islands, Spain

\bigskip

\end{center}

\date{ }

\bigskip

\begin{abstract}
We are proposing Tulczyjew's triple for 
contact dynamics. The most important ingredients of the triple, namely symplectic diffeomorphisms, special symplectic manifolds, and Morse families, are generalized to the contact framework. These geometries permit us to determine so-called generating family (obtained by merging a special contact manifold and a Morse family) for a Legendrian submanifold. Contact Hamiltonian and Lagrangian Dynamics are recast as Legendrian submanifolds of the tangent contact manifold. In this picture, the Legendre transformation is determined to be 
a passage between two different generators of the same Legendrian submanifold. A variant of contact Tulczyjew's triple  
is constructed for evolution  contact dynamics.
\textit{
}
\end{abstract}

\tableofcontents
\setlength{\parskip}{4mm}

\onehalfspacing

\section{Introduction}

Lagrangian Dynamics is generated by a Lagrangian function defined on the tangent bundle $T\mathcal{Q}$ of the configuration space of a physical system whereas Hamiltonian Dynamics is governed by a Hamiltonian function on the cotangent bundle $T^*\mathcal{Q}$ which is canonically symplectic, \cite{abraham1978foundations,holm2009geometric,Leon-book,LiMa87}. 
If a Lagrangian function is regular, that is, if it satisfies the Hessian condition, then the fiber derivative becomes a fibred local diffeomorphism from the tangent bundle to the cotangent bundle. In this case, the fiber derivative turns out be the Legendre transformation linking the Lagrangian and the Hamiltonian realizations of the physical system.

If a Lagrangian function happens to be degenerate then the fiber derivative fails to be a local diffeomorphism  since its  image space turns out only to be, in the best of the cases,  a proper submanifold of the cotangent bundle $T^*\mathcal{Q}$. That is, one only arrives at a presymplectic picture determined by some constraint functions. To deal with these constraints, Dirac proposed an algorithm, nowadays called Dirac-Bergmann algorithm, \cite{Di58,Di67}. This algorithm proposes a method to arrive at a submanifold (possibly smaller than the image space of the Legendre transformation) of the cotangent bundle where the Hamilton's equations becomes well-defined. In the final stage of the algorithm, one obtains the so-called Dirac bracket. 
There also exists a more geometric version of this approach  called the Gotay-Nester-Hinds algorithm \cite{GoNeHi78}. Inspiring from the tools introduced in \cite{GoNeHi78}, the Skinner-Rusk unified theory \cite{SkRu83} is establishing a unification of
Lagrangian and Hamiltonian formalisms on the Whitney sum of tangent and cotangent bundles. In this paper, we shall focus on the Tulczyjew approach for the Legendre transformations of singular Lagrangians.

\textbf{The Classical Tulczyjew's Triple.}
Tulczyjew's triple  is a commutative diagram linking three symplectic bundles namely $TT^* \mathcal{Q}$, $T^*T^* \mathcal{Q}$ and $T^*T\mathcal{Q}$ via symplectic diffeomorphisms, \cite{TuUr99}. This geometrization enables one to recast 
Lagrangian and Hamiltonian dynamical equations as Lagrangian submanifolds of the Tulczyjew symplectic space $TT^* \mathcal{Q}$ \cite{
tulczyjew1976soush}. Referring to this geometry 
the Legendre transformation is defined as a passage between two different generators of the same Lagrangian submanifold \cite{tulczyjew1972hamiltonian,Tu77}. This definition is free from the non-degeneracy requirement that is the Hessian condition. Evidently, this approach is in harmony with the creed by Weinstein "everything is a Lagrangian submanifold" 
\cite{weinstein1982symplectic}.  We reserve Section \ref{Sec-TT} for a brief summary of the Tulczyjew's triple and the Legendre transformation in this picture. 

Tulczyjew's triple is modified for many physical systems and it is carried to several geometric frameworks.  
For higher order Lagrangian dynamics, the triple is upgraded in \cite{deLeLa89,esen2018geometry}. For physical theories where the configuration space is a Lie group, the triple is examined in a series of works \cite{esen2014tulczyjew,esen2015tulczyjew,esen2021tulczyew,grabowska2016tulczyjew}.
The triple is examined for principal fiber bundles in \cite{GaGuMaMe,esen2020tulczyjew}. It is written for the vector bundle of $n$-vectors in 
\cite{GrGrUr14}. For the case of the field theories, we refer to an extensive but an incomplete list \cite{CaGuMa,cantrijn1999geometry,de2003tulczyjew,de1991special,SilviaBook,SilviaArt,echeverria2000geometry,
	grabowska2012tulczyjew, grabowska2013tulczyjew,roman2007k}, and for the higher order field theories see, for example, \cite{grabowska2015tulczyjew}.
	
\textbf{Contact Hamiltonian Dynamics.} A symplectic manifold must be even dimensional. An odd dimensional generalization of symplectic geometry is contact geometry \cite{arnold1989mathematical}. Hamiltonian Dynamics is available on this generalization as well, \cite{Br17,BrCrTa17,LeLa19,de2020review}. In the present work, as a generic model of contact manifold, we consider the extended cotangent bundle
\begin{equation}
\mathcal{T}^*\mathcal{Q}=T^*\mathcal{Q} \times \mathbb{R}
\end{equation}
with the contact one-form $\eta_\mathcal{Q}=dz-\theta_\mathcal{Q}$. Here, $z$ is the real variable, and $\theta_\mathcal{Q}$ is the pull-back of the canonical one-form on the cotangent bundle $T^*\mathcal{Q}$. In a coordinate free formulation, Contact Hamiltonian Dynamics generated by a Hamiltonian function $H$ is defined as 
	 \begin{equation}
\iota_{X^c_{H}}\eta_\mathcal{Q} =-H,\qquad \iota_{X^c_{H}}d\eta_\mathcal{Q} =dH-\mathcal{R}(H) \eta_\mathcal{Q},   \label{contact-intro}
\end{equation}%
where $\mathcal{R}$ is the Reeb field associated with the contact form $\eta_\mathcal{Q}$, and $\iota_{X^c_{H}}$ is the interior product. In Darboux' coordinates $(q^i,p_i,z)$, the contact Hamilton's equations are computed to be
\begin{equation}\label{conham-}
\dot{q}^i= \frac{\partial H}{\partial p_i}, \qquad \dot{p}_i = -\frac{\partial H}{\partial q^i}- 
p_i\frac{\partial H}{\partial z}, \quad \dot{z} = p_i\frac{\partial H}{\partial p_i} - H.
\end{equation}
Contact Hamiltonian Dynamics has different features from Classical Hamiltonian Dynamics. One interesting characteristic of Contact Hamiltonian Dynamics is the loss of conservation of Hamiltonian function along the motion \cite{gaset2020new,LeSa17}. Even, the canonical contact volume form is not preserved under the action of the contact Hamiltonian dynamics (see \cite{BrLeMaPa20}). The dissipative nature of Contact Hamiltonian Dynamics makes its proper for dissipative dynamical systems. We present here an incomplete list of some recent works along this direction: \cite{ciaglia2018contact,LeLa19,
LeSa17,gaset2020new}. Additionally, we can argue that contact framework looks proper for thermodynamics, see, for example, the following incomplete list:  \cite{Br17,bravetti2019contact,ghosh2019contact,
mrugala1991contact}.
Accordingly, we also wish to cite \cite{grmela2014contact,grmela2021multiscale,grmela1997dynamics,PavelkaKlikaGrmela2018} for the use of contact geometry for qualitative analysis of reversible-irreversible dynamics under GENERIC (an acronym for General Equation for Non-equilibrium
Reversible-Irreversible Coupling) formalism. 
 At the end of this paper, we shall provide some applications of our results on thermodynamics. 
We cite \cite{WaWaYa17,WaWaYa19} for some recent studies on variational aspects of Contact Hamiltonian Dynamics. 

There is evolution  contact Hamiltonian formalism on the extended cotangent bundle. In this case, for a Hamiltonian function $H$, evolution  contact Hamilton's equation is defined to be
\begin{equation}\label{evo-def--} 
\iota_{\varepsilon_H} \eta_\mathcal{Q}=0, \qquad \mathcal{L}_{\varepsilon_H}\eta_\mathcal{Q}=dH+\mathcal{R}(H)\eta_\mathcal{Q}.
\end{equation} 
Here, $\mathcal{L}_{\varepsilon_H}$ is the Lie derivative. According to the Cartan's formula, it is computed to be $\mathcal{L}_{\varepsilon_H}=d\iota_{\varepsilon_H}+\iota_{\varepsilon_H}d$. In Darboux' coordinates, the local picture of the evolution  contact Hamilton's equation \eqref{evo-def--} is computed to be 
\begin{equation}\label{evo-eq-}
	\dot{q}^i= \frac{\partial H}{\partial p_i}, \qquad \dot{p}_i = -\frac{\partial H}{\partial q^i}- 
	p_i\frac{\partial H}{\partial z}, \quad \dot{z} = p_i\frac{\partial H}{\partial p_i}.
\end{equation}
Note that if $0$ is a regular value of $H$, then the evolutionary vector field $\varepsilon_H$ is tangent to the hypersurface $H^{-1}(0)$ and, thus, to every Legendrian submanifold contained in $H^{-1}(0)$. Therefore, the integral curves of $\varepsilon_H$ may be interpreted as thermodynamical processes for a system with thermodynamical phase space the extended cotangent bundle. This fact was used in \cite{simoes2020contact} to discuss the relation between evolution contact dynamics and simple thermodynamical systems with friction.

\textbf{Contact Lagrangian Dynamics.} 
Corresponding Lagrangian formalism for Contact Hamiltonian Dynamics is also available in the literature under the name of Herglotz (or generalized Euler-Lagrange) formalism \cite{Herglotz}. See also \cite{de2020review}. In this paper, as the geometric framework, we consider the extended tangent bundle
\begin{equation}
\mathcal{T}\mathcal{Q}=T\mathcal{Q} \times \mathbb{R}.
\end{equation}
A Lagrangian function $L=L(q,\dot{q},z)$ on $\mathcal{T}\mathcal{Q}$ determines the Herglotz (generalized Euler-Lagrange) equations   as 
\begin{equation}\label{Herglotz-intro}
\frac{\partial L}{\partial q^i} - \frac{d}{dt}\Big(\frac{\partial L}{\partial {\dot q}^i} \Big)
+ \frac{\partial L}{\partial z}\frac{\partial L}{\partial {\dot q}^i} = 0. 
\end{equation}
 Here, $z$ is the real variable appearing in  $T\mathcal{Q} \times \mathbb{R}$. Evidently, if the Lagrangian function $L$ is independent of $z$ then Equations \eqref{Herglotz-intro} reduces to the classical Euler-Lagrange equations. 
If the Lagrangian function $L$ is non-degenerate that is if the rank of the Hessian matrix
 \begin{equation} \label{Leg-Trf-}
\left[ \frac{\partial ^2 L}{\partial \dot{q}^i \partial \dot{q}^i}\right]
\end{equation}
is maximum then the fiber derivative 
\begin{equation}\label{Herglotz-evo-1}
\mathbb{F}L^c: \mathcal
{T}\mathcal{Q} \longrightarrow  \mathcal
{T}^*\mathcal{Q}, \qquad (q^i,\dot{q}^i,z)\mapsto 
(q^i,\frac{\partial L}{\partial \dot{q}^i},z)
\end{equation}
turns out to be a local diffeomorphism. 
In this case, a direct calculation shows that the fiber derivative \eqref{Herglotz-evo-1} maps the Herglotz equations in \eqref{Herglotz-intro} to the contact  Hamilton's equations \eqref{conham-} once the Hamiltonian function is taken to be
\begin{equation}
H(q^i,p_i,z)=\dot{q}^ip_i-L(q,\dot{q},z).
\end{equation} 
If a Lagrangian function fails to be non-degenerate then, as in the symplectic case, the transformation \eqref{Herglotz-evo-1} fails to be a local isomorphism. Then the image space of $\mathbb{F}L^c$ can only be, in the best of the cases, a proper submanifold of the extended cotangent bundle. 
One way to deal with the constraints defining the image space is to employ a version of Dirac algorithm. This is recently studied in \cite{LeVa19b}. Another way is to generalize the unified formalism for contact dynamics as presented in a recent study \cite{LeGaLaRiRo20}.   
In this paper, our interest is to study the Legendre transformation for contact dynamics following the understanding of Tulczyjew.

\textbf{Goal of the Present Work.} The aim of this work is to define Legendre transformation between the Herglotz equations \eqref{Herglotz-intro} and the contact Hamilton's equations \eqref{conham-} by properly constructing a Tulczyjew's triple for the case of contact manifolds. We shall call this as contact Tulczyjew's triple. Such an attempt involves modifications of the ingredients of the classical Tulczyjew's triple to the contact geometry. Here, in the contact  case, the role of Tulczyjew's symplectic space $TT^*\mathcal{Q}$ will be played by the extended tangent bundle $\mathcal{T}\mathcal{T}^*\mathcal{Q}$ of the extended cotangent bundle. Accordingly,  Tulczyjew's triple 
will consist of the iterated extended bundles $\mathcal{T}^*\mathcal{T}\mathcal{Q}$, $\mathcal{T}\mathcal{T}^*\mathcal{Q}$, $\mathcal{T}^*\mathcal{T}^*\mathcal{Q}$ as well as contact transformation between them. Maybe, the most vital object in this picture is the introduction of the notion of special contact manifolds. This novel framework permits one to recast both the contact Lagrangian and the contact Hamiltonian dynamics as a Legendrian submanifold of $\mathcal{T}\mathcal{T}^*\mathcal{Q}$. This enables us to realize that Lagrangian and Hamiltonian functions as generating objects of the same Legendrian submanifold. So that, by merging special contact geometry with Morse family theory, the Legendre transformation for contact dynamics is defined to be a passage between two different generators of the same Legendrian submanifold. In addition, there is an (evolution) Hamiltonian flow on contact geometry preserving the energy but not the kernel of the contact form. There is also Lagrangian counterpart of this theory called evolution  Herglotz equations. 
By properly modifying contact  Tulczyjew's triple, the Legendre transformation for the evolution  Herglotz equations and the evolution  contact Hamilton's equations are obtained. In this theory, contact manifolds and Legendrian submanifolds are replaced by symplectic manifolds and Legendrian submanifolds, respectively. 
We call this geometry as the evolution  contact Tulczyjew's triple. 

\textbf{The content} of this work is as follows. The main body of the paper is consisting of three sections.  In Section \ref{Sec-TT},  for the sake of the completeness of the manuscript and in order to fix the notation, a brief summary of classical Tulczyjew's triple is given. Section \ref{Section-contact} is reserved for the basics on Contact Dynamics in both Hamiltonian and Lagrangian formulations. Section\ref{Section-TTc} is the one containing the novel results of the paper where the  Tulczyjew's triple is constructed for the contact and evolution  contact dynamics.

\section{The Classical Tulczyjew's triple}\label{Sec-TT}

     			\subsection{(Special) Symplectic Manifolds} \label{Ss-sss}~
		
		A manifold $\mathcal{P}$ is said to be symplectic if it is equipped with a non-degenerate closed two-form $\omega$ \cite{abraham1978foundations,daSilva-book,Leon-book}. In this case, $\omega$ is called a symplectic two-form. A diffeomorphism between two symplectic manifold is called a symplectic diffeomorphism if it respects the symplectic two-forms. 

\textbf{Submanifolds.}		Let $(\mathcal{P},\omega)$ be a symplectic manifold, and $\mathcal{S}$ be a submanifold of $\mathcal{P}$. We define the symplectic orthogonal complement
of $T\mathcal{S}$ as the vector subbundle of $T\mathcal{P}$ 
\begin{equation}
 T\mathcal{S}^{\bot}=\{X\in T\mathcal{P} ~: ~ \omega(X,Y)=0,~ \forall Y\in T\mathcal{S}\}.
\end{equation}
The rank of the tangent bundle $T\mathcal{P}$ is the sum of the ranks of the tangent bundle $T\mathcal{S}$ and its symplectic orthogonal complement $\ T\mathcal{S}^{\bot}$. We list some of the important cases. 
 \begin{itemize}
  \item $\mathcal{S}$ is called an isotropic submanifold if $T\mathcal{S}\subset T\mathcal{S}^{\bot}$. In this case, the dimension of $\mathcal{S}$ is less or equal to the half of the dimension of $\mathcal{P}$.
  \item $\mathcal{S}$ is called a coisotropic submanifold if $T\mathcal{S}^{\bot}\subset T\mathcal{S}$. In this case, the dimension of $\mathcal{S}$ is greater or equal to the half of the dimension of $\mathcal{P}$.
  \item $\mathcal{S}$ is called a Lagrangian submanifold if  $T\mathcal{S}=T\mathcal{S}^{\bot}$. In this case, the dimension of $\mathcal{S}$ is equal to the half of the dimension of $\mathcal{P}$. 
 \end{itemize}
Under a symplectic diffeomorphism, the image of a Lagrangian (isotropic, coisotropic) submanifold is a Lagrangian (resp. isotropic, coisotropic) submanifold.
		
\textbf{The Cotangent Bundle.}	 The generic examples of symplectic manifolds are cotangent bundles. To see this, consider a manifold $\mathcal{Q}$, and its cotangent bundle $T^{*}\mathcal{Q}$. The canonical (Liouville) one-form $\theta_\mathcal{Q}$ on $T^{*}\mathcal{Q}$ is defined, on a vector $X$ over $T^{*}\mathcal{Q}$  as
		\begin{equation} \label{can-Lio}
		\theta_\mathcal{Q}(X)= \langle \tau_{T^{*}\mathcal{Q}}(X), T\pi_\mathcal{Q}(X)\rangle.
		\end{equation}
	Here, $\tau_{T^{*}\mathcal{Q}}$ is the projection from the tangent bundle $TT^{*}\mathcal{Q}$ to its base manifold $T^{*}\mathcal{Q}$ whereas $T\pi_\mathcal{Q}$ is the tangent lift of the cotangent projection $\pi_{\mathcal{Q}}$. To be more precise, we present the following commutative diagram,
			\begin{equation}\label{one-form}
		\begin{tikzcd}
		&TT^*\mathcal{Q}\arrow[dr,"\tau_{T^*\mathcal{Q}}"]\arrow[dl,"T\pi_\mathcal{Q}",swap] \\
		T\mathcal{Q} \arrow[dr,"\tau_{\mathcal{Q}}",swap] && T^*\mathcal{Q}\arrow[dl,"\pi_{\mathcal{Q}}"] \\
		&\mathcal{Q}
		\end{tikzcd}
		\end{equation}
	Minus of the exterior derivative of  the canonical one-form $\theta_\mathcal{Q}$, that is, $\omega_\mathcal{Q}:=-d\theta_\mathcal{Q}$, is the canonical symplectic two-form on the cotangent bundle $T^{*}\mathcal{Q}$. 

\textbf{Hamiltonian Vector Fields.} Let $X$ be a vector field on the symplectic manifold $(T^*\mathcal{Q},\omega_\mathcal{Q})$. It is called a local Hamiltonian vector field if it preserves the symplectic two-form
\begin{equation}\label{X-Ham-pre}
\mathfrak{L}_X\omega_\mathcal{Q}=0.
\end{equation} 	
Poincar\'{e} Lemma assures us the existence of a local function $H$ satisfying the so-called Hamilton's equation
\begin{equation}\label{Ham-Eq}
\iota_{X_H}\omega_\mathcal{Q}=dH. 
\end{equation}
Notice that, in this realization, we denote the Hamiltonian vector field by $X_H$. If there exists a global Hamiltonian function $H$ then $X_H$ is called a (global) Hamiltonian vector field. 

If the dimension of $\mathcal{Q}$ is $n$ then the cotangent bundle $T^*\mathcal{Q}$ turns out to be a $2n$ dimensional manifold. In this case, the $n$-th power $\omega_\mathcal{Q}^n$ of the symplectic two-form is a non-vanishing top-form on $T^*\mathcal{Q}$. So that, it determines a volume form, the so-called symplectic volume. The identity \eqref{X-Ham-pre} gives that a Hamiltonian vector field preserves the symplectic volume. Further, the skew-symmetry of the symplectic two-form manifests the conservation of the Hamiltonian function $H$. If $H$ is the total energy then this is called the conservation of energy. 

A set of natural bundle coordinates $(q^i,p_i)$ on $T^{*}\mathcal{Q}$ are Darboux coordinates, meaning that the canonical one-form $\theta_\mathcal{Q}$ and the symplectic two-form $\omega_\mathcal{Q}$ are written as
		\begin{equation} \label{can-forms}
		\theta_\mathcal{Q}=p_idq^i, \qquad \omega_\mathcal{Q}=dq^i\wedge dp_i,
		\end{equation}
		respectively. In this realization, the Hamilton's equation \eqref{Ham-Eq} is computed to be
\begin{equation}\label{Ham-Eq-Loc}
\frac{dq^i}{dt}=\frac{\partial H}{\partial p_i},\qquad \frac{dp_i}{dt}=-\frac{\partial H}{\partial  q^i }.
\end{equation}

		\textbf{Special Symplectic Structures.} 
		Let $\mathcal{P}$ be a symplectic manifold carrying an exact symplectic two-form $\omega=-d\theta$, with $\theta$ is being a potential one-form. Assume also that, $\mathcal{P}$ is the total space of a fibre bundle $(\mathcal{P},\pi,\mathcal{Q})$. Here, $\mathcal{Q}$ is the base space,  and $\pi$ is the projection. A special symplectic structure is a quintuple 
		\begin{equation}  \label{sss}
		(\mathcal{P},\pi,\mathcal{Q},\theta,\phi)
		\end{equation} where $
		\phi$ is a fiber preserving symplectic
		diffeomorphism from $\mathcal{P}$ to the cotangent bundle $T^{\ast}\mathcal{Q}$ (see \cite{LaSnTu75, SnTu72}). 		
The symplectic diffeomorphism $\phi$ relating canonical symplectic two-form $\omega_{\mathcal{Q}}=-d\theta_{\mathcal{Q}}$ on $T^*\mathcal{Q}$ and the symplectic two-form $\omega=-d\theta$ satisfies both 
		\begin{equation}\label{symp-sss}
		\phi^*\theta_{\mathcal{Q}}=\theta,\qquad \phi^*\omega_{\mathcal{Q}}=\omega.
		\end{equation}
	Here, $\theta_{\mathcal{Q}}$ is the canonical one-form defined in \eqref{can-Lio}. 		
Accordingly, the symplectic 	diffeomorphism $\phi$ can be characterized by the following pairing
		\begin{equation} \label{chi-}
		\left\langle \phi(x),T\pi\circ X(x)\right\rangle = \left\langle
		\theta(x),X(x)\right\rangle
		\end{equation}
		for a vector field $X$ on $\mathcal{P}$, for any point $x$ in $\mathcal{P}$. To see this definition, simply evaluate $X$ with the pull-back one form $\phi^*\theta_{\mathcal{Q}}$ that is,
		\begin{equation}
		\begin{split}
		\langle \phi^*\theta_{\mathcal{Q}},X\rangle(x) &=\langle \theta_{\mathcal{Q}} (\phi (x)),T\phi \circ X (x)\rangle
		\\&= \langle \tau_{T^{*}\mathcal{Q}}\circ T\phi \circ X (x), T\pi_\mathcal{Q}\circ T\phi \circ X  (x) \rangle
		\\&=\langle \phi(x),T\pi\circ X(x)\rangle,
		\end{split}
		\end{equation}
		where we have employed the definition \eqref{can-Lio} of the canonical one-form in the second line whereas the identities $\tau_{T^{*}\mathcal{Q}}\circ T\phi \circ X =\phi $ and $\pi_\mathcal{Q}\circ \phi=\pi$ in the third line. 
	We exhibit a special symplectic manifold $  
		(\mathcal{P},\pi,\mathcal{Q},\theta,\phi)$ with the following commutative diagram	
		\begin{equation} \label{sss-}
		\xymatrix{ T^{\ast }\mathcal{Q}
			\ar[dr]_{\pi_{\mathcal{Q}}} &&\mathcal{P} \ar[ll]_{\phi}\ar[dl]^{\pi}
			\\&\mathcal{Q} }  
		\end{equation}
The two-tuple $(\mathcal{P},\omega)$ is called underlying symplectic manifold of the special symplectic structure. 
				
\subsection{Morse Families}\label{Sec-MF}

First of all, we will review the definition of a special kind of Lagrangian submanifolds of a cotangent bundle endowed with the canonical symplectic structure (for more details, see \cite{LiMa87,Marle-Jacobi}). 
Let $\mathcal{N}$ be a submanifold of a smooth manifold $\mathcal{Q}$ and $\Theta\vert_{\mathcal{N}}$ be a closed one-form defined on $\mathcal{N}$. Then we can consider a Lagrangian submanifold $\mathcal{S}_{\Theta\vert_{\mathcal{N}}}$ of $T^*\mathcal{Q}$ given by 
\begin{equation}
\mathcal{S}_{\Theta\vert_{\mathcal{N}}}=\{ \gamma\in T^*\mathcal{Q} ~ \vert ~ \pi_{\mathcal{Q}}(\gamma)\in \mathcal{N},   \big(T^*_{\pi_{\mathcal{Q}}(\gamma)}\big)\iota (\gamma)=\Theta\vert_{\mathcal{N}}(\pi_{\mathcal{Q}}(\gamma)) \} 
\end{equation}
where $\iota:\mathcal{N}\hookrightarrow \mathcal{Q}$ is the canonical inclusion. Note that for every point $q$ in $\mathcal{N}$ there exists an open subset $U\subseteq \mathcal{N}$ including  $q$, and a real valued smooth function $\Delta: U\mapsto \mathbb{R}$ such that $\Theta\vert_{\mathcal{N}}=d\Delta$. In this case, $\Delta$ is called the (local) generating function of $\mathcal{S}_{\Theta\vert_{\mathcal{N}}}$. 
Next, we will review the definition of Morse families and the Lagrangian submanifolds associated with them.

Consider a fiber bundle $(\mathcal{W},\tau,\mathcal{Q})$ where $\mathcal{W}$ is the total space 
of dimension $n+K$ whereas $\mathcal{Q}$ is the base manifold of dimension $n$. Here, $\tau$ is the bundle projection. The vertical bundle $V\tau$ over the manifold $\mathcal{W}$ is a vector subbundle of $T\mathcal{W}$ containing vectors that belong to the kernel of the tangent mapping $T\tau$ that is,
\begin{equation}
V\tau=\{X\in T\mathcal{W}~: ~ T\tau \circ X =0\}.
\end{equation}
The conormal bundle $V^0\tau$ is the space of covectors in $T^*\mathcal{W}$ annihilating vectors in the vertical bundle $V\tau$. 
A real-valued function $E$ on the total space of a fiber bundle $(\mathcal{W},\tau,\mathcal{Q})$ is called Morse family (or  generating family) \cite{LiMa87,SnTu72} if 
\begin{equation}
T _\gamma {\rm im}(dE)+T _\gamma V^0\tau=T _\gamma T^*\mathcal{W},
\end{equation} 
for all $\gamma$ in the intersection ${\rm im}(dE)\cap V^0\tau$. 

\textbf{Generating Lagrangian Submanifolds.} A Morse family $E$  defined on $(\mathcal{W},\tau,\mathcal{Q})$ generates an immersed  Lagrangian submanifold of the cotangent bundle $T^*\mathcal{Q}$ as
\begin{equation}\label{LagSub}
\mathcal{S}=\big\{z \in T^{\ast}\mathcal{Q}~:~T^{\ast}\tau(z)=dE ( w ), ~\forall w\in \mathcal{W}, \tau(w)=\pi_{\mathcal{Q}}(z)
\big \}. 
\end{equation} 
The inverse of this statement is also true. So we state the following theorem which is generalizing the well-known Poincar\'{e} Lemma for non-horizontal Lagrangian submanifolds. Accordingly, it is called as generalized Poincaré Lemma or Maslov-Hörmander Theorem \cite{Be11,LiMa87,We77}.
\begin{theorem}
For a Lagrangian submanifold of a cotangent bundle $T^*\mathcal{Q}$, there always exists, at least locally, a Morse  family  $E$ generating $\mathcal{S}$. 
\end{theorem}
We picture the Lagrangian submanifold $\mathcal{S}$ generated by a Morse family $E$ on $(\mathcal{W},\tau,\mathcal{Q})$ as follows.
\begin{equation} \label{Morse-pre}
 \xymatrix{
 	\mathbb{R}& \mathcal{W} \ar[d]^{\tau}\ar[l]^{E}& T^*\mathcal{Q}\ar@(ul,ur)^{\mathcal{S}}   \ar[d]^{\pi_{\mathcal{Q}}}\\ &
 	\mathcal{Q}  \ar@{=}[r]& \mathcal{Q}
 }
 \end{equation}
Given a Lagrangian submanifold, its Morse family generator is far from being unique. For example, one may find a Morse family with less number of fiber variables generating the same Lagrangian submanifold. This procedure is called reduction of Morse family. See \cite{Be11} for further discussions on this subject. 

We merge a Morse family $E$ defined on $(\mathcal{W},\tau,\mathcal{Q})$ and a special symplectic structure $  
		(\mathcal{P},\pi,\mathcal{Q},\theta,\phi)$ depicted diagrammaticality in \eqref{sss-}. This permits us to define a Lagrangian submanifold of the symplectic manifold $(\mathcal{P},\omega)$  \cite{TuUr99}. Here is the diagram,
		\begin{equation} \label{Morse-pre-sss}
 \xymatrix{
 	\mathbb{R}& \mathcal{W} \ar[d]^{\tau}\ar[l]^{E}& T^*\mathcal{Q}\ar@(ul,ur)^{\mathcal{S}}   \ar[dr]_{\pi_{\mathcal{Q}}}
 	&&
 	\mathcal{P} \ar[ll]_{\phi}\ar[dl]^{\pi}  \ar@(ul,ur)^{ \mathcal{S}_{E}}
 	\\ &
 	\mathcal{Q}  \ar@{=}[rr]&& \mathcal{Q}
 }
 \end{equation}
Here, $\mathcal{S}$ is the Lagrangian submanifold of $T^*\mathcal{Q}$ given in \eqref{LagSub}, and the inverse of the symplectic diffeomorphism $\phi$ maps $\mathcal{S}$ to a Lagrangian submanifold $\mathcal{S}_{E}$ of $\mathcal{P}$. That is, $\phi(\mathcal{S}_{E})=\mathcal{S}$. 	

\textbf{Local Picture.} Now we present local realization of these discussions. 
Let $(q^i)$ be local coordinates on $\mathcal{Q}$ and we consider the induced local coordinates $(q^i,\epsilon^a)$ on the the total space $\mathcal{W}$. In this image, a function $E$ is called a Morse family if the rank of the following matrix 
\begin{equation}
\left( \frac{\partial ^{2}E}{\partial q^i \partial \epsilon^a}, \qquad \frac{%
	\partial ^{2}E}{\partial \epsilon^a \partial \epsilon^b}\right)  \label{MorseCon} 
\end{equation}%
is equal to $K$. In this case, the Lagrangian submanifold \eqref{LagSub} generated by $E$ is viewed locally as 
\begin{equation} \label{MFGen}
\mathcal{S}=\left \{\Big(q^i,\frac{\partial E}{\partial q^j}(x,\epsilon)\Big )\in T^*\mathcal{Q}: \frac{\partial E}{\partial \epsilon^a}(x,\epsilon)=0\right \}.
\end{equation} 
Note that the dimension of the submanifold $\mathcal{S}$ is half of the dimension of the cotangent bundle $T^*\mathcal{Q}$ and the canonical symplectic two-form $\omega_\mathcal{Q}$ vanishes on $\mathcal{S}$. 

\subsection{Merging Two Special Symplectic Structures}~
\label{Merge-Sec}

Let $(\mathcal{P},\omega)$ be an exact symplectic manifold. Assume that $\mathcal{P}$ admits two different fiber bundle structures denoted by $(\mathcal{P},\pi,\mathcal{Q})$ and $(\mathcal{P},\pi',\mathcal{Q}')$. Further, let these fibrations lead to two special symplectic symplectic structures denoted by $(\mathcal{P},\pi,\mathcal{Q},\theta,\phi)$ and $(\mathcal{P},\pi',\mathcal{Q}',\theta',\phi')$, respectively. We merge these two special symplectic structures in one diagram as follows:
\begin{equation}\label{TT-P}
\begin{tikzcd}
T^*\mathcal{Q}\arrow[rdd,"\pi_{\mathcal{Q}}"]&&\mathcal{P} \arrow[ll,"\phi",swap]\arrow[rr,"\phi'"]\arrow[rdd,"\pi'"]\arrow[ldd,"\pi",swap]  && T^*\mathcal{Q}'\arrow[ldd,"\pi_{\mathcal{Q}'}",swap] \\\\
& \mathcal{Q} && \mathcal{Q}'
\end{tikzcd}
\end{equation}
This is the most abstract realization of classical Tulczyjew's triple.  Notice that, in this geometry, the symplectic two-form $\omega$ on $\mathcal{P}$ admits two different potential one forms $\theta'$ and $\theta$ so that
 \begin{equation}\label{oto}
\omega=-d\theta'=-d\theta. 
\end{equation}
Further, by employing the canonical symplectic two-forms $\omega_{\mathcal{Q}}=-d\theta_{\mathcal{Q}}$ and $\omega_{\mathcal{Q}'}=-d\theta_{\mathcal{Q}'}$ over the corresponding cotangent bundles $T^*\mathcal{Q}$ and $T^*\mathcal{Q}'$, the following properties hold
\begin{equation}
\phi^*(\theta_{\mathcal{Q}})= \theta,\qquad \phi^*(\omega_{\mathcal{Q}})= \omega,\qquad
\phi'^*(\theta_{\mathcal{Q}'})= \theta',\qquad \phi'^*(\omega_{\mathcal{Q}'})= \omega.
\end{equation}

Consider the product manifold $\mathcal{Q}\times \mathcal{Q}'$ and let $\pi$ and $\pi'$ be the maps defined in~\eqref{TT-P}. Assume that the image space of 
\begin{equation}\label{chi}
\chi:\mathcal{P}\longrightarrow  \mathcal{Q}\times \mathcal{Q}',\qquad p\mapsto (\pi(p),\pi'(p))  
\end{equation}
is an embedded submanifold ${\rm im}\chi=\mathcal{N}$ of $\mathcal{Q}\times \mathcal{Q}'$. Further, we assume that the map $\chi:\mathcal{P}\mapsto \mathcal{N}$ is a surjective submersion. 
The equalities \eqref{oto} give  that the difference $\theta'-\theta$ is a closed one-form. Using Poincar\'{e} Lemma, we have that for every point $p$ in $\mathcal{P}$ there exists an open subset $U_p\subseteq \mathcal{P}$ containing $p$  and a smooth function $\Delta_p$ on $U_p$ such that  
\begin{equation}\label{Delta--}
d\Delta_p=(\theta'-\theta)\big\vert_{U_p}.
\end{equation} 
 The vertical bundle with respect to the fibration $\chi$ is precisely the intersection of the vertical bundles with respect to $\pi$ and $\pi'$ that is, $V\chi=V\pi\cap V\pi'$. Notice that $\theta'-\theta$ is a $\chi$-basic one-form. In fact, due to the fiber preserving character of the symplectic diffeomorphisms $\phi$ and $\phi'$, the potential one-forms $\theta$ and $\theta'$ are taking values in the conormal bundles $V^0\pi$ and  $V^0\pi'$, respectively. 
We deduce that the difference of potential one-forms takes values in the conormal bundle $V^0\chi$ since
\begin{equation}
\theta'-\theta\in \Gamma(V^0\pi') + \Gamma(V^0\pi)=
\Gamma(V^0\pi'+ V^0\pi)=\Gamma(V \pi' \cap V\pi)^0=\Gamma(V^0\chi),
\end{equation}
where $\chi$ is the fibration in \eqref{chi}. 
So, there exists a unique closed one-form $\Theta\vert_{\mathcal{N}}$  such that $\chi^*\Theta\vert_{\mathcal{N}}=\theta'-\theta$. In fact, $\Delta_p$ given in \eqref{Delta--} is  a $\chi$-basic function and, therefore, there exists a smooth function $\Delta$ on $\chi(U_p)$ satisfying   $d\Delta=\theta\big\vert_{U_p}$. The following diagram illustrates the above situation.
\begin{equation}\label{TT-P--}
\begin{tikzcd}
U_p \subseteq \mathcal{P} \arrow[rd,"\Delta_p"]
\arrow[dd,"\chi",swap] \\
& \mathbb{R}
\\
\chi(U_p) \arrow[ru,"\Delta",swap]
\end{tikzcd} 
\end{equation}
Next, we consider the following Lagrangian submanifold of the canonical  symplectic manifold $(T^*(\mathcal{Q}'\times \mathcal{Q}),\omega_{\mathcal{Q}'\times \mathcal{Q}})$ given by
\begin{equation}\label{bar-S}
\bar{\mathcal{S}}=\big\{ \bar{\gamma}\in T^*(\mathcal{Q}'\times \mathcal{Q}): \pi_{\mathcal{Q}'\times \mathcal{Q}}(\bar{\gamma})=(q',q)\in \mathcal{N},\quad \bar{\gamma}\big\vert_{T_{(q',q)}\mathcal{N}}=\Theta\vert_{\mathcal{N}}(q',q)\big
\}.
\end{equation}
Consider the symplectic diffeomorphism 
\begin{equation}\label{Psi}
\Psi:\big(T^*(\mathcal{Q}'\times \mathcal{Q}),\omega_{\mathcal{Q}'\times \mathcal{Q}}\big)\longrightarrow \big(T^*\mathcal{Q}'\times T^* \mathcal{Q}, \omega_{\mathcal{Q}'}\ominus\omega_{\mathcal{Q}} \big),\qquad (\mu_{q'},\mu_{q})\mapsto (\mu_{q'},-\mu_{q})
\end{equation}
where $\omega_{\mathcal{Q}'\times \mathcal{Q}}$ is the canonical form on the cotangent bundle whereas  
\begin{equation*}
	\omega_{\mathcal{Q}'}\ominus\omega_{\mathcal{Q}} ={\rm pr}_2^*(\omega_{\mathcal{Q}'})-{\rm pr}_1^*(\omega_{\mathcal{Q}}).
\end{equation*}
Referring to \eqref{Psi}, we map the Lagrangian submanifold $\bar{\mathcal{S}}$ in \eqref{bar-S} to a Lagrangian submanifold $\mathcal{S}=\Psi(\bar{\mathcal{S}})$ of $T^*\mathcal{Q}'\times T^* \mathcal{Q}$. Then we will prove the following result.
\begin{theorem}
The ${\rm graph}(\phi\circ (\phi')^{-1})$ of the symplectic diffeomorphism 
\begin{equation}
\phi\circ (\phi')^{-1}:T^*\mathcal{Q}'\mapsto T^*\mathcal{Q}
\end{equation}
pictured in (\ref{TT-P}) is an open subset of the Lagrangian submanifold $\mathcal{S}=\Psi(\bar{\mathcal{S}})$. 
\end{theorem}

\textbf{Proof.} Evidently, ${\rm graph}(\phi\circ (\phi')^{-1})$ is a Lagrangian submanifold of $T^*\mathcal{Q}'\times T^* \mathcal{Q}$ so that it is enough to show that ${\rm graph}(\phi\circ (\phi')^{-1})$ is a subset of $\mathcal{S}$. In fact, we will see that 
\begin{equation*}
{\rm graph}(\phi\circ (\phi')^{-1})=\big\{(\mu_{q'},-\mu_q)\in T^*\mathcal{Q}'\times T^*\mathcal{Q} : q=\pi((\phi')^{-1}(\mu_{q'})), (\mu_{q'},-\mu_q)\big\vert_{T_{(q',q)}\mathcal{N}}=\Theta\vert_{\mathcal{N}}(q',q)\big\}
\end{equation*}
which implies the result. Suppose that $(\mu_{q'},-\mu_q)$ be an element of ${\rm graph}(\phi\circ (\phi')^{-1})$, then we have that $\phi\circ (\phi')^{-1}(\mu_{q'})$ is precisely $\mu_q$. Assume that $(\phi')^{-1}(\mu_{q'})=p$ in $\mathcal{P}$ so that $\phi(p)=-\mu_q$. This gives 
\begin{equation}
\pi_{\mathcal{Q}'}(\phi(p))=q \, \Rightarrow \, \pi(p)=q \, \Rightarrow \, \pi\circ (\phi')^{-1}(\mu_{q'})=q.
\end{equation}
Next, we will see that $(\mu_{q'},-\mu_q)\big\vert_{T_{(q',q)}\mathcal{N}}=\Theta\vert_{\mathcal{N}}(q',q)$. Using \eqref{chi-}, it follows that 
\begin{equation}
\begin{split}
\phi'(p)&=\mu_{q'}\, \Rightarrow \, \theta'(p)=T^*_p\pi'(\mu_{q'})
\\
\phi(p)&=-\mu_{q}\, \Rightarrow \, -\theta(p)=T^*_p\pi(\mu_{q})
\end{split}
\end{equation}
therefore one has that
\begin{equation}
\theta'(p)-\theta(p)= T^*_p\chi\circ T^*_{(q',q)}\iota(\mu_{q'},\mu_{q}).
\end{equation}
But since 
\begin{equation}
T^*_p\chi \circ \Theta\vert_{\mathcal{N}}(q',q)=\theta'(p)-\theta(p)
\end{equation}
we conclude that 
\begin{equation}
\begin{split}
T^*_p\chi \circ \Theta\vert_{\mathcal{N}}(q',q))&=T^*_p\chi\circ  T^*_{(q',q)}\iota(\mu_{q'},\mu_{q}) ,
\\
\Theta\vert_{\mathcal{N}}(q',q)&=T^*_{(y',y)}\iota(\mu_{q'},\mu_{q})
=(\mu_{q'},\mu_q)\big\vert_{T_{(q',q)}\mathcal{N}}.
\end{split}
\end{equation}
Conversely, suppose that $\gamma=(\mu_{q'},-\mu_q)$ in $T^*\mathcal{Q}'\times T^*\mathcal{Q}$ 
 and 
$q=\pi((\phi')^{-1}(\mu_{q'})) $, $(\mu_{q'},-\mu_q)\big\vert_{T_{(q',q)}\mathcal{N}}=\Theta\vert_{\mathcal{N}}(q',q)$. If $p=(\phi')^{-1}(\mu_{q'})$, we will see that $\phi(p)=-\mu_q$ which implies that $(\mu_{q'},-\mu_{q})$ in ${\rm graph}(\phi\circ (\phi')^{-1})$. We have that $\phi'(p)=\mu_{q'}$ which gives $\theta'(p)=T^*_p\pi'(\mu_{q'})$. On the other hand, 
 \begin{equation}
 \begin{split}
T^*_{(q',q)}\iota(\mu_{q'},\mu_{q})&=\Theta\vert_{\mathcal{N}}(q',q)
\\
T^*_{p}\chi \circ T^*_{(q',q)}\iota(\mu_{q'},\mu_{q})&=\theta'(p)-\theta(p)\\
T^*(\iota\circ \chi)(\mu_{q'},\mu_{q})&=\theta'(p)-\theta(p)
\\
T^*_p\pi'(\mu_{q'})+T^*_p\pi(\mu_{q})=\theta'(p)-\theta(p).
\end{split}
 \end{equation}
 Thus, as $T^*_p\pi'(\mu_{q'})=\theta'(p)$, we conclude that 
 \begin{equation}
 T^*_p\pi(-\mu_{q})=\theta(p)\, \Rightarrow \, \pi(p)=-\mu_{q}.
 \end{equation}

To illustrate the Legendre transformation in this geometry we present the following discussion.
Assume the existence of the triple in \eqref{TT-P}. Start with the left wing of the triple by considering a real valued function $F$ on $\mathcal{Q}$. This determines a Lagrangian submanifold ${\rm im}(dF)$ of $T^*\mathcal{Q}$. By employing the inverse of the symplectic diffeomorphism $\phi$ to ${\rm im}(dF)$, one arrives at a Lagrangian submanifold 
$\mathcal{S}_{F}$ of the symplectic manifold $\mathcal{P}$. This submanifold can also be defined in terms of the potential one-form $\theta$ as 
\begin{equation*}
\mathcal{S}_{F}=\{p\in \mathcal{P}~:~ d(F \circ \pi)(p)=\theta(p)\}.
\end{equation*}
The Legendre transformation in terms of Tulczyjew \cite{Tu77,TuUr99} is to determine a generating family for $\mathcal{S}_{\mathcal{F}}$ referring to the right wing of the triple \eqref{TT-P}. To have that, first map $\mathcal{S}_{F}$ to a Lagrangian submanifold of $T^*\mathcal{Q}'$ by means of the symplectic diffeomorphism $\phi'$. As a manifestation of the generalized Poincar\'{e} lemma, there exits a Morse family $E'$ on a fiber bundle $(\mathcal{W},\tau,\mathcal{Q})$ generating $\phi(\mathcal{S}_{F})$. 

\subsection{The Classical Tulczyjew's triple} \label{CTT}

In this subsection, we draw Tulczyjew's triple for Classical Dynamics assuming a configuration manifold $\mathcal{Q}$. This is to construct the triple \eqref{TT-P} by replacing $\mathcal{Q'}$ with the tangent bundle $T\mathcal{Q}$ whereas replacing $\mathcal{Q}$ with the cotangent bundle  $T^*\mathcal{Q}$. On the upper level of \eqref{TT-P}, this results with the iterated bundles $T^*T\mathcal{Q}$, $TT^*\mathcal{Q}$ and $T^*T^*\mathcal{Q}$ in order. See that, being cotangent bundles both $T^*T\mathcal{Q}$ and $T^*T^*\mathcal{Q}$ are symplectic. We now establish the symplectic structure on $TT^*\mathcal{Q}$ admitting two potential one forms.

\noindent \textbf{Symplectic structure on $TT^*\mathcal{Q}$.} Consider the canonical symplectic manifold $T^{\ast}Q$ equipped with the exact symplectic two-form $\omega_\mathcal{Q}=-d\theta_\mathcal{Q}$. The derivation $i_{T}$ takes the symplectic two-form $\omega_\mathcal{Q}$ on $T^*\mathcal{Q}$ to a one-form on $TT^*\mathcal{Q}$ as
\begin{equation*}
  i_{T}\omega _{\mathcal{Q}}(Y)=\omega _{\mathcal{Q}}(\tau_{TT^*\mathcal{Q}} (Y),T\tau_{T^{\ast}Q}(Y)),
\end{equation*}
for any tangent vector $Y$ on $TT^*\mathcal{Q}$. 
Here, $\tau_{TT^*\mathcal{Q}}$ is the tangent bundle projection $TTT^*\mathcal{Q}$ to $TT^*\mathcal{Q}$ whereas $T\tau_{T^{\ast}\mathcal{Q}}$ is the tangent mapping of the bundle projection $\tau_{T^{\ast}\mathcal{Q}}:TT^*\mathcal{Q}\mapsto T^*\mathcal{Q}$. We define two one-forms on $TT^*\mathcal{Q}$ as  
		\begin{equation} \label{tet2}
		\vartheta _{1}=-i_{T}\omega_{\mathcal{Q}}, \qquad  \vartheta _{2}= d_{T}\theta _{\mathcal{Q}}= i_{T}d\theta_{Q}+di_{T}\theta _{\mathcal{Q}}. 
		\end{equation}
where the derivation $d_{T}$ is the commutator $i_{T}d+di_{T}$. The exterior derivatives of these one-forms results with a symplectic two-form on $TT^*\mathcal{Q}$ defined to be 
\begin{equation}\label{d-T-Omega-Q}
d_T\omega_\mathcal{Q}=-d\vartheta_1=-d\vartheta_2.
\end{equation}
We record this in the following theorem \cite{Tu77}.
\begin{theorem}
The tangent bundle $TT^*\mathcal{Q}$ is a symplectic manifold with the symplectic two-form $d_T\omega_\mathcal{Q}$ in \eqref{d-T-Omega-Q} admitting two potential one-forms given in \eqref{tet2}.
\end{theorem}

Consider the Darboux' coordinates $(q^i,p_i)$ on the cotangent bundle $T^{\ast}\mathcal{Q}$. 
In terms of the induced local coordinate chart $(q^i,p_i;\dot{q}^i,\dot{p}_i)$ on the tangent bundle $TT^*\mathcal{Q}$, the potential one-forms in \eqref{tet2} are computed to be
\begin{equation}
\vartheta_{1}= - i_{T}\omega _{\mathcal{Q}}=\dot{p}_idq^i-\dot{q}^i dp_i, \qquad \vartheta_{2}= d_{T}\theta _{\mathcal{Q}}= \dot{p}_i dq^i + p_i d\dot{q}^i  \label{thets}.
\end{equation}
Notice that, in this case, the symplectic two-form turns out to be 
\begin{equation}
	d_T\omega_\mathcal{Q}=dq^i\wedge d\dot{p}_i+d\dot{q}^i\wedge dp_i.
\end{equation}
Note that, the value $\vartheta_{2}-\vartheta_{1}$ is an exact one-form. Actually, it is the exterior derivative of coupling function  $i_T\theta_{\mathcal{Q}}:TT^*\mathcal{Q}\mapsto \mathbb{R}$. 

One can arrive at the tangent bundle symplectic  two-form $d_T\omega_\mathcal{Q}$ on $TT^*\mathcal{Q}$ as the complete lift of the canonical symplectic two-form $\omega_\mathcal{Q}$ on $T^*\mathcal{Q}$. To have this, from \cite{Leon-book,YaIs73} we recall the definition of the complete lift of a differential form. The complete lift of the canonical one-form $\theta_\mathcal{Q}=p_i dq^i$ is computed to be 
 \begin{equation}
 \theta_\mathcal{Q}^C=\dot{p}_i dq^i+p_i d\dot{q}^i.
 \end{equation}
Complete lift of forms commutes with exterior
derivative. So that, we compute
 \begin{equation}
 \omega_\mathcal{Q}^C=(-d \theta_\mathcal{Q})^C=dq^i\wedge d\dot{p}_i+d\dot{q}^i\wedge dp_i=d_T\omega_\mathcal{Q}.
 \end{equation}
As manifested in the display, we conclude that $\omega_\mathcal{Q}^C=d_T\omega_\mathcal{Q}$. 

By recalling Diagram \ref{one-form}, notice that  $TT^*\mathcal{Q}$ admits two bundle structures. It is a vector bundle over $T\mathcal{Q}$ with respect to the vector bundle projection $T\pi_{\mathcal{Q}}$, and it is a vector bundle over $T^*\mathcal{Q}$ with vector bundle projection $\tau_{T^*\mathcal{Q}}$. So, we can consider the projection
\begin{equation}
\sigma:TT^*\mathcal{Q}\longrightarrow \mathfrak{N}=T\mathcal{Q}\times_{\mathcal{Q}} T^*\mathcal{Q},\qquad Z\mapsto (T\pi_{\mathcal{Q}}(Z),\tau_{T^*\mathcal{Q}}(Z)) 
\end{equation}
from $TT^*\mathcal{Q}$ to the Whitney sum $ \mathfrak{N}=T\mathcal{Q}\times_{\mathcal{Q}} T^*\mathcal{Q}$. The function $i_T\theta_{\mathcal{Q}}$ is $\sigma$-basic and it induces a smooth function $\Delta$ on  $ \mathfrak{N}$ which is just the coupling function
\begin{equation}\label{Delta}
\Delta=\dot{q}^i p_i.
\end{equation} 
In addition, $ \mathfrak{N}$ is a submanifold of the product manifold $ T\mathcal{Q}\times T^*\mathcal{Q}$. So we can consider the Lagrangian submanifold $\mathcal{S}_{\theta _{\mathfrak{N}}}$ induced by the exact one-form $\theta _{\mathfrak{N}}=d\Delta$ as in Subsection \ref{Sec-MF}. Moreover, using the construction presented in Subsection \ref{Merge-Sec}, we deduce that $\mathcal{S}_{\theta _{\mathfrak{N}}}$ is the graph of a symplectic diffeomorphism $\psi$ between $T^*T\mathcal{Q}$ and $T^*T^*\mathcal{Q}$. In terms of the Darboux' coordinates $(q^i,\dot{q}^i,a_i, \dot{a}_i)$ on $T^*T\mathcal{Q}$, it is computed to be
\begin{equation} \label{shuttle-loc}
\psi: T^*T\mathcal{Q}\longrightarrow T^*T^*\mathcal{Q}, \qquad  (q^i,\dot{q}^i,a_i, \dot{a}_i)  \mapsto
          (q^i,\dot{a}_i,a_i,-\dot{q}^i).
\end{equation}

\textbf{The canonical involution on $TT\mathcal{Q}$.} The iterated tangent bundle $TT\mathcal{Q}$ is a tangent bundle with the base manifold $T\mathcal{Q}$ along with the tangent bundle projection $\tau_{T\mathcal{Q}}$. It is possible to show that $TT\mathcal{Q}$ can be written as a vector bundle over $T\mathcal{Q}$ apart from the canonical tangent bundle fibration. This fibration is achieved by the tangent mapping $T\tau_\mathcal{Q}$ of the projection $\tau_\mathcal{Q}$. 
		To manifest this, we plot here the double vector bundle structure of the iterated tangent bundle $TT\mathcal{Q}$ in the following commutative diagram
		\begin{equation}\label{tangent-rhombic}
		\begin{tikzcd}
		&TT\mathcal{Q}\arrow[dl,"\tau_{T\mathcal{Q}}",swap]\arrow[dr,"T\tau_\mathcal{Q}"] \\
		T\mathcal{Q} \arrow[dr,"\tau_{\mathcal{Q}}",swap] && T\mathcal{Q}\arrow[dl,"\tau_{\mathcal{Q}}"] \\
		&\mathcal{Q}
		\end{tikzcd}
		\end{equation}
		Referring to this double bundle structure \cite{abraham1978foundations,Godbillon,Ur96}, consider a differential mapping $\Gamma=\Gamma(t,s)$ from an open domain of $\mathbb{R}^2$ to $\mathcal{Q}$. The differential of $\Gamma$ with respect to $t$, at $t=0$, results with a curve $\gamma(s)$ lying in the tangent bundle $T\mathcal{Q}$ depending on the free variable $s$. If one further takes the differential of the curve $\gamma(s)\subset T\mathcal{Q}$ with respect to $s$, at $s=0$, then arrives at a vector in $TT\mathcal{Q}$. Accordingly, the canonical involution $\kappa_\mathcal{Q}$ on $TT\mathcal{Q}$ is defined by changing the order of differentiations
\begin{equation} \label{involution}
\kappa:TT\mathcal{Q}\longrightarrow TT\mathcal{Q}: \left. \frac{d}{dt} \right \vert_{t=0} \left. \frac{d}{ds} \right \vert_{s=0} \Gamma(t,s) \longrightarrow \left.  \frac{d} {ds} \right \vert_{s=0} \left. \frac{d} {dt} \right \vert_{t=0}\Gamma(t,s).
\end{equation} 
An observation gives that the involution $\kappa $ is changing the order of the fibrations 
		\begin{equation} \label{can-inv}
			\tau_{T\mathcal{Q}}\circ \kappa =T\tau_\mathcal{Q}, \qquad T\tau_\mathcal{Q} \circ \kappa =\tau_{T\mathcal{Q}}.
		\end{equation}

		\textbf{Pairing between $TT^*\mathcal{Q}$ and $TT\mathcal{Q}$.}	We now establish a pairing between an element $Z$ in $TT^{\ast }\mathcal{Q}$ and an element     $W$ in $TT\mathcal{Q}$ such that $T\pi_{\mathcal{Q}}(Z) = T\tau_{\mathcal{Q}}(W)$ and that $\pi_{\mathcal{Q}}\circ z =\tau_{\mathcal{Q}}\circ  v$. Recall that, there is a curve $z(t)$ in $T^{\ast }\mathcal{Q}$ so that $Z=\dot{z}(0)$ and, there is a curve $v(t)$ in $T\mathcal{Q}$ so that $W=\dot{v}(0)$. In this framework, the pairing is defined as
		\begin{equation} \label{pairing-tilde}
		\langle \bullet,\bullet \rangle^{\widetilde{}} :TT^\ast \mathcal{Q} \times TT\mathcal{Q}\longrightarrow 
		\mathbb{R}, \qquad \langle Z,W\rangle^{\widetilde{}} =\left .\frac{d}{dt}
		\langle z(t),v(t) \rangle \right \vert
		_{t=0}.
		\end{equation}
Here, the pairing on the right hand side is the one between $T^*\mathcal{Q}$ and $T\mathcal{Q}$. 
		
		In the induced coordinates $(q^i,\dot{q}^i,{q}^{\prime i},\dot{q}^{\prime i})$ on the iterated tangent bundle $TT\mathcal{Q}$, the fibration in \eqref{tangent-rhombic} read
		\begin{equation}\label{local-tau-TQ-Q}
		\tau_{T\mathcal{Q}}(q^i,\dot{q}^i,{q}^{\prime i},\dot{q}^{\prime i})=(q^i,\dot{q}^i), \qquad T\tau_{\mathcal{Q}}(q^i,\dot{q}^i,{q}^{\prime i},\dot{q}^{\prime i})=(q^i, {q}^{\prime i}),
		\end{equation}
		whereas the canonical involution $\kappa $ in \eqref{can-inv} is computed to be 
		\begin{equation*}
		\kappa (q^i,\dot{q}^i,{q}^{\prime i},\dot{q}^{\prime i})=(q^i,{q}^{\prime i},\dot{q}^i,\dot{q}^{\prime i}).
		\end{equation*}
		Coordinate expression of the pairing \eqref{pairing-tilde} is as follows. Let us choose coordinates on $TT^*\mathcal{Q}$ as $Z=(q^i,p_i,\dot{q}^i,\dot{p}_i)$, and coordinates $W=(q^i,{q}^{\prime i},\dot{q}^i,\dot{q}^{\prime i})$ on $TT\mathcal{Q}$ then, 
		\begin{equation*}
		\langle Z,W\rangle^{\widetilde{}}=\big\langle (q^i,p_i,\dot{q}^i,\dot{p}_i),(q^i,{q}^{\prime i},\dot{q}^i,\dot{q}^{\prime i})\big\rangle^{\widetilde{}} = p_i\dot{q}^{\prime i}+{q}^{\prime i}\dot{p}_i.
		\end{equation*} 
		
Now we define two special symplectic structures for the symplectic manifold $(TT^*\mathcal{Q},d_T\omega_\mathcal{Q})$. One is to $T^*T\mathcal{Q}$ and the other is to $T^*T^*\mathcal{Q}$. We assume that, being cotangent bundles, $T^*T\mathcal{Q}$ and $T^*T^*\mathcal{Q}$ are equipped with the canonical symplectic forms $\omega_{T\mathcal{Q}}=-d\theta_{T\mathcal{Q}}$ and $\omega_{T^*\mathcal{Q}}=-d\theta_{T^*\mathcal{Q}}$, respectively.

		\textbf{Left Wing of the Tulczyjew's triple.} Start with defining the vector fibration morphism 
		\begin{equation}\label{alpha}
		\alpha : TT^*\mathcal{Q} \longrightarrow T^*T\mathcal{Q}, \qquad  \langle \alpha (Z) ,W \rangle
		=\langle Z, \kappa(W) \rangle ^{\widetilde{}},  
		\end{equation}
		where $\kappa$ is the canonical involution defined in \eqref{can-inv} whereas the pairing on the right hand side is the one given in \eqref{pairing-tilde}. Here, the pairing of the left hand side is the canonical pairing between $T^*T\mathcal{Q}$ and $TT\mathcal{Q}$. Locally, one has that
		\begin{equation}\label{alpha-locally}
	\alpha (q^i,p_i,\dot{q}^i,\dot{p}_i)=(q^i,\dot{q}^i, \dot{p}_i, p_i).		
			\end{equation}
See that, $\alpha$ is a symplectic diffeomorphism by satisfying $\alpha ^{\ast }\omega _{T\mathcal{Q}}=d_T\omega_\mathcal{Q}$. Here, $\omega _{T\mathcal{Q}}$ is the canonical symplectic two-form on $T^*T\mathcal{Q}$, and $d_T\omega_\mathcal{Q}$ is the lifted symplectic two-form \eqref{d-T-Omega-Q} on $TT^*\mathcal{Q}$. So that we arrive at a special symplectic structure
				\begin{equation}\label{sss-1}
		(TT^{\ast }Q,T\pi _{\mathcal{Q}},T\mathcal{Q},\vartheta
		_{2},\alpha ) 
			\end{equation}		
		where $\vartheta
		_{2}$ is the one-form in \eqref{thets}, and $T\pi _{\mathcal{Q}}$ is the tangent lift of the cotangent bundle projection $\pi _{\mathcal{Q}}$. We include this in the following diagram for future reference
		\begin{equation}\label{ssp-diagram-left}
\begin{tikzcd}
T^*T\mathcal{Q} \arrow[dr,"\pi_{T\mathcal{Q}}",swap]&& TT^*\mathcal{Q}\arrow[ll,"\alpha",swap]\arrow[dl,"T\pi_\mathcal{Q}"]
\\
& T\mathcal{Q} &
\end{tikzcd}
\end{equation}

	\textbf{Right Wing of the Tulczyjew's triple.}  The nondegeneracy of the canonical symplectic two-form $\omega_{\mathcal{Q}}$ leads to the existence of a (musical) diffeomorphism
		\begin{equation} \label{beta}
		\beta: TT^*\mathcal{Q}\longrightarrow T^*T^*\mathcal{Q}, \qquad Z\mapsto -\omega^\flat_\mathcal{Q}(Z)=-\omega_\mathcal{Q}(Z,\bullet).
	\end{equation}
Locally, one has that
	\begin{equation}
							\beta(q^i,p_i,\dot{q}^i,\dot{p}_i)=(q^i,p_i,\dot{p}_i,-\dot{q}^i).
						\end{equation}
Note that, $\beta$ is a symplectomorphism by satisfying $\beta^{\ast }\omega _{T^*\mathcal{Q}}=d_T\omega_\mathcal{Q}$. Here, $\omega _{T^*\mathcal{Q}}$ is the canonical symplectic two-form on $T^*T^*\mathcal{Q}$. So that we arrive at the following special symplectic structure
						\begin{equation}\label{sss-2}
		(TT^{\ast }Q,\tau _{T^{\ast }Q},T^{\ast }Q
		,\vartheta _{1},\beta).
		\end{equation}
	Here,  $\vartheta _{1}$ is the potential one-form in \eqref{tet2} and that $\tau _{T^{\ast }Q}$ is the tangent bundle projection. The diagram is for the future reference.
		\begin{equation}\label{ssp-diagram-right}
\begin{tikzcd}
TT^*\mathcal{Q} \arrow[dr,"\tau_{T^*\mathcal{Q}}",swap]\arrow[rr,"\beta"]&& T^*T^*\mathcal{Q}\arrow[dl,"\pi_{T^*\mathcal{Q}}"]
\\
& T^*\mathcal{Q} &
\end{tikzcd}
\end{equation}	
 The Tulczyjew's triple for Classical Mechanics is the following commutative diagram merging the two special symplectic structures \eqref{ssp-diagram-left} and \eqref{ssp-diagram-right} in one picture 
		\begin{equation} \label{TT}
		\xymatrix{T^{\ast }T\mathcal{Q} \ar[dr]_{\pi_{T\mathcal{Q}}}&&TT^{\ast
			}Q\ar[dl]^{T\pi_{\mathcal{Q}}}
			\ar[rr]^{\beta} \ar[dr]_{\tau_{T^{\ast
					}Q}} \ar[ll]_{\alpha}&&T^{\ast }T^{\ast
			}Q\ar[dl]^{\pi _{T^{\ast}Q}}
			\\&T\mathcal{Q}\ar[dr]_{\tau_{\mathcal{Q}}}
			&&T^{\ast}Q\ar[dl]^{\pi_{\mathcal{Q}}} \\&&\mathcal{Q}}
		\end{equation}
		where $\beta$ and $\alpha$ are 
symplectomorphisms such that 
\begin{equation*}
\beta^*\omega_{T^*\mathcal{Q}}=\omega^C_\mathcal{Q}=-d\vartheta_1, \qquad  \alpha ^*\omega_{T\mathcal{Q}}=\omega^C_\mathcal{Q}=-d\vartheta_2.
\end{equation*}
By employing the composition of two symplectic diffeomorphisms $\alpha^{-1}$ and $\beta$, we arrive at 
a symplectic diffeomorphism $\psi=\beta\circ \alpha^{-1}$. This is exactly the one in \eqref{shuttle-loc}.

\subsection{The Legendre Transformation}

Let us now realized Euler-Lagrange equations generated by a Lagrangian function $L$ on the tangent $T\mathcal{Q}$ as a Lagrangian submanifold of $TT^*\mathcal{Q}$ using the left side of Tulczyjew's triple (\ref{TT}). The image space of the exterior derivative $dL$ is a Lagrangian submanifold of the manifold $T^*T\mathcal{Q}$. By employing the inverse of $\alpha$, we transfer this Lagrangian submanifold to a Lagrangian submanifold $\mathcal{S}_L$ of the manifold $TT^*\mathcal{Q}$. 
In terms of the local coordinates $(q^i,p_i,\dot{q}^i,\dot{p}_i)$ on $TT^*\mathcal{Q}$, $\mathcal{S}$ is computed to be
\begin{equation} \label{LagSubEL}
\mathcal{S}_L=\alpha ({\rm im}dL)=\left\{\left(q^i ,\frac{\partial L}{\partial \dot{q}^i}, \dot{q}^i ,\frac{\partial L}{\partial  q^i}\right)\in TT^*\mathcal{Q}:L=L(q^i,\dot{q}^i)  \right\}\subset TT^*\mathcal{Q}.
\end{equation}
Dynamics determined by the Lagrangian submanifold $\mathcal{S}_L$ is computed simply by taking the time derivatives of the base coordinates $\big(q^i ,{\partial L}/{\partial \dot{q}^i}\big)$
 and equate them to the fiber coordinates given by $(\dot{q}^i ,{\partial L}/{\partial  q^i})$, respectively. By this, we arrive at the Euler-Lagrange equations
\begin{equation}
\frac{dq^i}{dt}=\dot{q}^i,\qquad  \frac{d}{dt}\Big( 
\frac{\partial L}{\partial \dot{q}^i}
\Big)=\frac{\partial L}{\partial  q^i}.
\end{equation}
Notice that, for this realization, we have not asked any regularity conditions for the Lagrangian function. That is, this theory is valid for singular Lagrangians as well.

\textbf{The Legendre Transformation of Euler-Lagrange Equations.}	We now 	generate the Lagrangian submanifold $\mathcal{S}_L$ given in \eqref{LagSubEL} referring to  the right wing of the Tulczyjew's triple (\ref{TT}). Let us remark once more that, this is the Legendre transformation in the understanding of Tulczyjew. To have that, at first consider the Morse family over the Pontryagin bundle as
\begin{equation} \label{energy-1}
E:T\mathcal{Q} \times_{\mathcal{Q}} T^*\mathcal{Q}\longrightarrow \mathbb{R},\qquad (u,\varsigma)\mapsto \langle \varsigma, u \rangle -L (u) 
\end{equation}
where $\langle \bullet , \bullet \rangle$ is the canonical pairing between $T^*\mathcal{Q}$ and $T\mathcal{Q}$. 
In coordinates, we write the Morse family as 
\begin{equation}
E(q^i,\dot{q}^i,p_i)= p_i \dot{q}^i-L(q,\dot{q}).
\end{equation} 
We remark that the matrix 
\begin{equation}
\left( \frac{\partial ^{2}E}{\partial \dot{q}^i \partial q^j}, \qquad \frac{\partial ^{2}E}{\partial \dot{q}^i \partial p_j
 },
\qquad
\frac{
	\partial ^{2}E}{\partial  \dot{q}^i \partial \dot{q}^j }\right) 
\end{equation}
has maximal rank. For the choice of sign conventions, we consider the minus of the Morse family, that is, $-E$. So, according to the definition in \eqref{MFGen}, the minus of the Morse family $-E$ generates a Lagrangian submanifold of the cotangent bundle $T^*T^*\mathcal{Q}$ as 
\begin{equation} \label{LagSubEL-T*T*}
\mathcal{S}=\left\{\left(q^i ,p_i, -\frac{\partial E}{\partial  q^i }, -\frac{\partial E}{\partial p_i} \right)\in T^*T^*\mathcal{Q}: \frac{\partial E}{\partial  \dot{q}^i}=0  \right\}\subset T^*T^*\mathcal{Q}.
\end{equation} 
In order to transfer this Lagrangian submanifold to $TT^*\mathcal{Q}$, we merge the Pontryagin bundle with the right wing of the Tulczyjew's triple (\ref{TT}) as follows 
\begin{equation}\label{Ham-Morse-Gen-}
\xymatrix{TT^{\ast }\mathcal{Q} \ar[rr]^{ \beta} \ar[dr]_{\tau_{T^{\ast }\mathcal{Q}}}
&&T^{\ast }T^{\ast }\mathcal{Q}\ar[dl]^{\pi _{T^{\ast}\mathcal{Q}}} &T\mathcal{Q} \times_{\mathcal{Q}} T^*\mathcal{Q} \ar[d]^{{\rm pr}_2}\ar[rr]^{\qquad -E} &&  \mathbb{R}
\\&T^{\ast}\mathcal{Q}  \ar@{=}[rr]&& T^{\ast}\mathcal{Q}}
\end{equation}
This permits us to arrive at the following Lagrangian submanifold of $TT^*\mathcal{Q}$
\begin{equation} \label{LagSubEL-}
\mathcal{S}_{-E}=\beta^{-1}(\mathcal{S})=\left\{\left(q^i ,p_i, \frac{\partial E}{\partial p_i}, -\frac{\partial E}{\partial  q^i }\right)\in TT^*\mathcal{Q}: \frac{\partial E}{\partial  \dot{q}^i}=0  \right\}\subset TT^*\mathcal{Q}.
\end{equation}
A direct computation proves that the Lagrangian submanifold $\mathcal{S}_{-E}$ in (\ref{LagSubEL-}) and the Lagrangian submanifold $\mathcal{S}_L$ in (\ref{LagSubEL}) are the same. So that the Legendre transformation is achieved. 
If the Lagrangian function is non-degenerate then from the equation
\begin{equation}
\frac{\partial E}{\partial \dot{q}^i}(q,\dot{q},p)=p_i-\frac{\partial L}{\partial\dot{q}^i}(q,\dot{q})=0
\end{equation} 
one can explicitly determine the velocity $\dot{q}^i$ in terms of the momenta $(q^i,p_i)$. In other words, for a non-degenerate Lagrangian function $L=L( q ,\dot{q})$ the fiber derivative
\begin{equation}
\mathbb{F}L:T\mathcal{Q}\longrightarrow T^*\mathcal{Q},\qquad  (q^i ,\dot{q}^j) \longrightarrow \big( q^i ,\frac{\partial L}{\partial \dot{q}^i}( q ,\dot{q})\big)
\end{equation}
is a local diffeomorphism. In this case, the Morse family $E$ can be reduced to a well-defined Hamiltonian function
\begin{equation} \label{canHam}
H(q^i,p_i)=p_i~\dot{q}^i( q,p)  -L\big( q ,\dot{q}( q ,p)\big)
\end{equation}
on $T^*\mathcal{Q}$. 

\textbf{Inverse Legendre Transformation.} The inverse Legendre transformation is also possible in a similar way. This time, one starts with a Hamiltonian system $(T^*\mathcal{Q},\omega_{\mathcal{Q}},H)$ where $H$ is a Hamiltonian function. See that, in this notation Hamiltonian vector field $X_H$  defined in \eqref{Ham-Eq} is determined through
				\begin{equation}\label{bemol-H}
\beta \circ X_H =-dH.
\end{equation} 
Notice that the Lagrangian submanifold determined by the equality \eqref{bemol-H} is written in coordinates as
\begin{equation} \label{LagSubEL-2}
\mathcal{S}_{-H}=\left\{\left(q^i ,p_i, \frac{\partial H}{\partial p_i}, -\frac{\partial H}{\partial  q^i }\right)\in TT^*\mathcal{Q}\right\}\subset TT^*\mathcal{Q}.
\end{equation}
Evidently, this Lagrangian submanifold is precisely determining the Hamilton's equation \eqref{Ham-Eq-Loc}. In the present picture, the inverse Legendre transformation is to generate the Lagrangian submanifold \eqref{LagSubEL-2} by referring to the right wing of the triple.  	
If the Hamiltonian function is not regular then one needs to employ a Morse family
\begin{equation}
F
:T\mathcal{Q} \times_{\mathcal{Q}} T^*\mathcal{Q}\longrightarrow \mathbb{R},\qquad (u,\varsigma)\mapsto \langle \varsigma, u \rangle -H  (\varsigma) .
\end{equation}
So, if we consider the Pontryagin bundle over $T\mathcal{Q}$ and we proceed as in the previous subsection, we will obtain the inverse Legendre transformation. 
\section{Contact Dynamics} \label{Section-contact}

\subsection{Contact Manifolds}

A  $(2n+1)-$dimensional manifold $\mathcal{M}$ is called contact manifold if it is equipped with a contact one-form $\eta$ satisfying $d\eta^n
\wedge \eta \neq 0$,  \cite{arnold1989mathematical,LiMa87}. We denote a contact manifold by a two-tuple $(\mathcal{M},\eta)$. The Reeb vector field $\mathcal{R}$ is the unique vector field  satisfying 
\begin{equation}
\iota_{\mathcal{R}}\eta =1,\qquad \iota_{\mathcal{R}}d\eta =0.
\end{equation}
At each point of the manifold $\mathcal{M}$, the kernel of the contact form $\eta$ determines the contact structure $H\mathcal{M}$. The complement of this structure, denoted by $V\mathcal{M}$,  is determined by the kernel of the exact two-form $d\eta$. These give the following decomposition of the tangent bundle 
\begin{equation}\label{decomp-TM}
T\mathcal{M}=H\mathcal{M}\oplus V\mathcal{M}, \qquad H\mathcal{M}=\ker \eta, ~ V\mathcal{M}=\ker d\eta.
\end{equation}
Here, $H\mathcal{M}$ is a vector subbundle of rank $2n$.  
The restriction of $d\eta$ to $H\mathcal{M}$ is non-degenerate so that  $(H\mathcal{M}, d\eta)$ is a symplectic vector bundle over $\mathcal{M}$. The rank of $V\mathcal{M}$ is $1$ and it is generated by  the Reeb field $\mathcal{R}$. 

\textbf{Contactization.}
It is possible to arrive at a contact manifold starting from a symplectic manifold. To have this, consider a symplectic manifold $\mathcal{P}$ admitting an integer symplectic two form $\omega$. Introduce the principal circle (quantization) bundle 
\begin{equation}\label{contactization}
S^{1}\rightsquigarrow (\mathcal{M},\eta )\overset{{\rm pr} }{\longrightarrow }%
(\mathcal{P},\omega). 
\end{equation}%
The contact one-form  on $\mathcal{M}$ 
is the connection one-form associated with a principal connection on the principal $S^{1}$-bundle $pr:\mathcal{M}\mapsto \mathcal{P}$ with curvature $\omega$. This procedure is called contactization. 

Another example of a contact manifold can be obtained from an exact symplectic manifold as follows. 
Consider a trivial line bundle over a manifold given by $\mathcal{Q}\times \mathbb{R}\mapsto \mathcal{Q}$. The first jet bundle, denoted by $\mathcal
{T}^*\mathcal{Q}$ is diffeomorphic to the product space $T^*\mathcal{Q}\times \mathbb{R}$ that is,
\begin{equation} \label{ext-cot}
\mathcal{T}^*\mathcal{Q}=T^*\mathcal{Q}\times \mathbb{R}.
\end{equation}
We call this space as the extended cotangent bundle.
There exist two projections
\begin{equation}\label{Sec-Jet}
\begin{split}
\pi^1_\mathcal{Q}&:\mathcal{T}^*\mathcal{Q}=T^*\mathcal{Q}\times \mathbb{R}\longrightarrow T^*\mathcal{Q},\qquad (\zeta,z)\mapsto \zeta
\\
\pi^0_\mathcal{Q}&:\mathcal{T}^*\mathcal{Q}=T^*\mathcal{Q}\times \mathbb{R}\longrightarrow \mathcal{Q}, \qquad (\zeta,z)\mapsto \pi_\mathcal{Q}(\zeta),
\end{split}
\end{equation}
where $\pi_\mathcal{Q}$ is the cotangent bundle projection whereas $z$ is the standard coordinate on $\mathbb{R}$. Referring to the fibration defined by $\pi^1_\mathcal{Q}$, we have the following globally trivial contactization of the canonical symplectic manifold
\begin{equation}\label{contactization-T*}
\mathbb{R}\rightsquigarrow (\mathcal{T}^*\mathcal{Q},\eta_{\mathcal{Q}} )\overset{\pi^1_\mathcal{Q}}{\longrightarrow }
(T^*\mathcal{Q},\omega_{\mathcal{Q}}).
\end{equation}
Here, the contact one-form on the jet bundle $\mathcal{T}^*\mathcal{Q}$ is defined to be
\begin{equation}\label{eta-Q}
\eta_{\mathcal{Q}}:=dz-\theta_{\mathcal{Q}}
\end{equation}
where $\theta_{\mathcal{Q}}$ is the canonical one-form \eqref{can-Lio} on the cotangent bundle $T^*\mathcal{Q}$. Notice that, we have employed abuse of notation by identifying $z$ and $\theta_{\mathcal{Q}}$ with their pull-backs on the total space $\mathcal{T}^*\mathcal{Q}$. The previous construction also works if we replace $T^*\mathcal{Q}$ by an arbitrary exact symplectic manifold $\mathcal{P}$ and, in such a case, we obtain a contact structure on the product manifold $\mathcal{P}\times \mathbb{R}$.  
There exist Darboux' coordinates $(q^i,p_i,z)$ on $\mathcal{T}^*\mathcal{Q}$, where $i$ is running from $1$ to $n$. In these coordinates, the contact one-form and the Reeb vector field are computed to be
\begin{equation}
\eta _{\mathcal{Q}}= d z -  p_i d q^i, \qquad \mathcal{R}=\frac{\partial}{\partial z},
\end{equation}
respectively. Notice that, in this realization, the horizontal bundle is generated by the vector fields
\begin{equation}\label{Horizontal-space}
H\mathcal{T}^*\mathcal{Q}=span\{\xi_i,\xi^i\},\qquad \xi_i=\frac{\partial}{\partial q^i} + p_i \frac{\partial}{\partial z},~\xi^i=\frac{\partial}{\partial p_i}.
\end{equation}
It is important to note that these generators are not closed under the Jacobi-Lie bracket that is, 
\begin{equation}
[\xi^{i},\xi_{j}]=\delta^i_j\mathcal{R},
\end{equation}
where $\delta^i_j$ stands for the Kronecker delta. The Darboux' theorem manifests that local picture presented in this subsection is generic for all contact manifolds of dimension $2n+1$.

\textbf{Musical Mappings.} For a contact manifold $(\mathcal{M}, \eta)$, there is a musical isomorphism $\flat$ from the tangent bundle $T\mathcal{M}$ to the cotangent bundle $T^*\mathcal{M}$ defined to be 
\begin{equation}\label{flat-map}
\flat:T\mathcal{M}\longrightarrow T^*\mathcal{M},\qquad v\mapsto \iota_vd\eta+\eta(v)\eta.
\end{equation} 
This mapping takes the Reeb field $\mathcal{R}$ to the contact one-form $\eta$. We denote the inverse of this mapping by $\sharp$. Referring to this, we define a bivector field $\Lambda$ on  $\mathcal{M}$ as
\begin{equation}\label{Lambda}
\Lambda(\alpha,\beta)=-d\eta(\sharp\alpha, \sharp \beta). 
\end{equation}
The couple $(\Lambda,-\mathcal{R})$ induces a Jacobi structure \cite{Kirillov-Local-Lie,Lichnerowicz-Jacobi,Marle-Jacobi}. This is a manifestation of the equalities 
\begin{equation}
[\Lambda,\Lambda]=-2\mathcal{R}\wedge \Lambda, \qquad [\mathcal{R},\Lambda]=0,
\end{equation}
where the bracket is the Schouten–Nijenhuis bracket. We cite \cite{BrGrGr17,LeMaPa97,LiMa87,Lichnerowicz-Jacobi} for more details on the Jacobi structure associated with a contact one-form. Referring to the bivector field $\Lambda$ we introduce the following musical mapping 
\begin{equation}\label{Sharp-Delta}
\sharp_\Lambda: T^*\mathcal{M}\longrightarrow T\mathcal{M}, \qquad  \alpha\mapsto \Lambda(\alpha,\bullet)= \sharp \alpha - \alpha(\mathcal{R})  \mathcal{R}. 
\end{equation}
Evidently, the mapping $\sharp_\Lambda$ fails to be an isomorphism. Notice that, the kernel is spanned by the contact one-form $\eta$. So that, the image space of $\sharp_\Lambda$ is precisely the horizontal bundle $H\mathcal{M}$ exhibited in \eqref{Horizontal-space}. 

In terms of the Darboux coordinates $(q^i,p_i,z)$, we compute the image of a one-form in $T^*\mathcal{M}$ by $\sharp_\Lambda$ as
\begin{equation}\label{zap}
\sharp_\Lambda:\alpha_i dq^i + \alpha^i dp^i + 
u dz\mapsto 
\alpha^i \frac{\partial}{\partial q^i}-(\alpha_i + p_i u)\frac{\partial}{\partial p_i}  + \alpha^ip_i  \frac{\partial}{\partial z}. 
\end{equation}

\textbf{Symplectization.} 
The symplectization of a contact manifold $(\mathcal{M},\eta)$ is the
symplectic manifold $(\mathcal{M} \times \mathbb{R},d(e^t\eta))$, where $t$ denotes the standard coordinate on $\mathbb{R}$
factor. In this case, $\mathcal{M}\times\mathbb{R}$ is said to be the symplectification of $\mathcal{M}$.  
The inverse of this assertion is also true. That is, if $(\mathcal{M} \times \mathbb{R},d(e^t\eta))$ is a symplectic manifold, then $(\mathcal{M},\eta)$  turns out to be contact.

\subsection{Submanifolds of Contact Manifolds}\label{sec-sub}

Let $(\mathcal{M},\eta)$ be a contact manifold. Recall the associated bivector field $\Lambda$ defined in \eqref{Lambda}. Consider a linear subbundle $ \Xi$ of the tangent bundle $T\mathcal{M}$ (that is, a distribution on $\mathcal{M}$). We define the contact complement of $\Xi$ as
\begin{equation}
\Xi^\perp : = \sharp_ \Lambda(\Xi^o),
\end{equation}
 where the sharp map on the right hand side is the one in   \eqref{Sharp-Delta} and   $\Xi^o$ is the annihilator of  $\Xi$.
Let $\mathcal{N}$ be a submanifold of $\mathcal{M}$. We say that $\mathcal{N}$ is:
    \begin{itemize}
        \item \emph{Isotropic} if $T\mathcal{N}\subseteq {T\mathcal{N}}^{\perp }$.
        \item \emph{Coisotropic} if $T\mathcal{N}\supseteq {T\mathcal{N}}^{\perp }$.
        \item \emph{Legendrian} if $T\mathcal{N}= {T\mathcal{N}}^{\perp }$.
    \end{itemize}

Assume that a submanifold $\mathcal{N}$ of a contact  manifold $\mathcal{M}$ is defined to be the zero level set of $k$ real smooth functions 
$\phi_a:U\to \mathbb{R}$. We determine $k$ vector fields $Z_a = \sharp_\Lambda (d \phi_a)$. The image space of these vector fields are spanning the contact complement
\begin{equation}
       {T\mathcal{N}}^{\perp} = {\rm span} \{ Z_a \; | \; a=1, \dots, k \}. 
\end{equation}
In this geometry, $\mathcal{N}$ is coisotropic if and only if, $Z_a(\phi_b)=0$ for all  $a,b$. According to the local computation in \eqref{zap}, we have that $\mathcal{N}$ is coisotropic if and only if
\begin{equation}\label{za2}
\frac{\partial \phi_a}{\partial p_i}\big( \frac{\partial \phi_b}{\partial q^i} + p_i \frac{\partial \phi_b}{\partial z} \big)
-(\frac{\partial \phi_a}{\partial q^i} + p_i \frac{\partial \phi_a}{\partial z}) \frac{\partial \phi_b}{\partial p_i}
  = 0.
\end{equation}
Referring to this local observation, one can easily prove that a submanifold $\mathcal{N}$ of $\mathcal{M}$ is Legendrian if and only if it is a maximal integral manifold of $\ker \eta$. In this case, the dimension of $\mathcal{N}$ must be $n$ (see \cite{LeLa19,deLeon2020infinitesimal,LiMa87}).

 
\textbf{Generating Functions.} Consider the first order jet bundle  $\mathcal{T}^*\mathcal{Q}$ endowed with the contact structure given in  \eqref{eta-Q}. 
Let $F$ be a real valued function on the base manifold $\mathcal{Q}$. Its first prolongation is a section of the bundle $\pi^0_\mathcal{Q}$ displayed in  \eqref{Sec-Jet} that is,
\begin{equation}\label{j1F}
\mathcal{T}^* F:\mathcal{Q}\longrightarrow \mathcal{T}^*\mathcal{Q}=T^*\mathcal{Q}\times \mathbb{R},\qquad q\mapsto (dF(q),F(q)).
\end{equation}
The image space of the first prolongation $\mathcal{T}^* F$ is  a Legendrian  submanifold of $\mathcal{T}^*\mathcal{Q}$. The converse of this assertion is also true, that is, if the image space of a section $\sigma$ of $\pi^0_\mathcal{Q}$ is a Legendrian submanifold then it is  the first prolongation of a function $F$.  Evidently, this is not the only way to obtain a Legendrian submanifold. 

Consider, for example, a Morse family $E$ defined on a smooth bundle $(\mathcal{W},\tau,\mathcal{Q})$ according to Subsection \ref{Sec-MF}. Then, referring to the definition \eqref{LagSub}, we define a Lagrangian submanifold $\mathcal{S}$ of the cotangent bundle $T^*\mathcal{Q}$. In the light of the first jet prolongation in \eqref{j1F}, we lift this Lagrangian submanifold to a Legendrian submanifold of the contact manifold $\mathcal{T}^*\mathcal{Q}$. To see this, consider a local system of coordinates $(q^i)$ on the base manifold $\mathcal{Q}$, and the induced coordinates  $(q^i,\epsilon^a)$ on the total space $\mathcal{W}$. Then, referring to the Darboux' coordinates on $\mathcal{T}^*\mathcal{Q}$, the Legendrian submanifold $\mathcal{N}$  generated by a Morse family $E=E(q,\epsilon)$ is computed to be 
\begin{equation} \label{MFGen-C}
\mathcal{N}  =\left \{\Big(q^i,\frac{\partial E}{\partial q^i}(q,\epsilon),E(q,\epsilon)\Big )\in \mathcal{T}^*\mathcal{Q}: \frac{\partial E}{\partial \epsilon^a}=0\right \}\subset \mathcal{T}^*\mathcal{Q}.
\end{equation} 

On the other hand, the lift of a Legendrian submanifold to the symplectification is a Lagrangian submanifold. In fact, $\mathcal{N}$ is a Legendrian submanifold of a contact manifold $(\mathcal{M},\eta)$ if and only if $\mathcal{N}\times \mathbb{R}$ is a
Lagrangian submanifold of the symplectic manifold $(\mathcal{M}\times \mathbb{R},d(e^t\eta))$, see \cite{IbLeMaMa97}.

\subsection{Contact Diffeomorphisms and Contact  Hamiltonian systems}

Let $(\mathcal{M}_1,\eta_1)$ and $(\mathcal{M}_2,\eta_2)$ be two contact manifolds. A diffeomorphism $\varphi$ from $\mathcal{M}_1$ to $\mathcal{M}_2$ is said to be a contact diffeomorphism (or contactomorphism) if it preserves the contact structures that is, $T\varphi(\ker \eta_1)=\ker \eta_2$. In terms of the contact forms, a contact diffeomorphism $\varphi$ is the one satisfying 
   \begin{equation}\label{Cont-Dif}
        \varphi^*\eta_2 = \mu \eta_1.
    \end{equation}
 where $\mu $ is a non-zero conformal factor. 
To manifest the existence of this conformal factor, a mapping $\varphi$ satisfying \eqref{Cont-Dif} is also called as conformal  contact diffeomorphism. In this understanding, the contact mapping is denoted by a two-tuple $(\varphi,\mu)$.

For a contact manifold $(\mathcal{M},\eta)$, we  denote the
group of contact diffeomorphisms \cite{Ba97} by 
\begin{equation}\label{Diff-con}
{\rm Diff}_{con} ( \mathcal{M}) =\left\{ \varphi \in {\rm Diff} ( \mathcal{%
M} ) :\varphi ^{\ast }\eta =\mu  \eta, \quad \mu  \in \mathcal{F}%
( \mathcal{M}) \right\} .
\end{equation}
Here, ${\rm Diff} ( \mathcal{
M})$ is standing for the group of all diffeomorphism on $ \mathcal{%
M}$.
A vector field on the contact
manifold $\left( \mathcal{M},\eta \right) $ is a contact vector field (called also as infinitesimal conformal contactomorphism) if
it generates one-parameter group of contact diffeomorphisms. We have that the space of contact vector fields is given by
\begin{equation}\label{algcon}
\mathfrak{X}_{con} ( \mathcal{M} ) =\left\{ X\in \mathfrak{X} ( 
\mathcal{M}) :\mathcal{L}_{X}\eta  =-\lambda \eta ,\quad \lambda
\in \mathcal{F} ( \mathcal{M} ) \right\} .  
\end{equation}
Sometimes a contact vector field is denoted by a two-tuple $(X,\lambda)$ in order to manifest the existence of the conformal factor $\lambda$. 
In order to discuss the geometry of $\lambda $, we perform the following observation. This permits us to introduce Hamiltonian dynamics on the present framework as well (for more details, see \cite{Br17,BrCrTa17,de2021hamilton,LeLa19}).

\textbf{(Contact) Hamiltonian Vector Fields.}
For a real valued function $H$ on a contact manifold $(\mathcal{M},\eta)$, there corresponds a contact vector field $X^c_H$ defined as follows
\begin{equation}
\iota_{X^c_{H}}\eta =-H,\qquad \iota_{X^c_{H}}d\eta =dH-\mathcal{R}(H) \eta,   \label{contact}
\end{equation}%
where $\mathcal{R}$ is the Reeb vector field. Here, $H$ is called the (contact) Hamiltonian function and $X^c_H$ is called the (contact)  Hamiltonian vector field. We denote a contact Hamiltonian system as a three-tuple $(\mathcal{M},\eta,H)$ where $(\mathcal{M},\eta)$ is a contact manifold and $H$ is a smooth real function on $M$. A direct computation determines the conformal factor for a given Hamiltonian vector fields as
\begin{equation}\label{L-X-eta}
\mathcal{L}_{X^c_{H}}\eta =
d\iota_{X^c_{H}}\eta+\iota_{X^c_{H}}d\eta= -\mathcal{R}(H)\eta.
\end{equation}
That is, $\lambda=\mathcal{R}(H)$. 

In this realization, the contact Jacobi bracket of two smooth functions on $\mathcal{M}
$ is defined by
\begin{equation}\label{cont-bracket}
\{F,H\}^c=\iota_{[X^c_F,X^c_H]}\eta, 
\end{equation}
 where $X_F$ and $X_H$ are Hamiltonian vectors fields determined through \eqref{contact}.  Here, $%
\left[ \bullet,\bullet \right]$ is the Lie bracket of vector
fields. Then, the identity 
\begin{equation} 
 -\left[
X^c_{K},X^c_{H}\right]=X^c_{\left\{ K,H\right\}^c}  
\end{equation}
 establishes the
isomorphism 
\begin{equation}
\left( \mathfrak{X}_{con}\left( \mathcal{M}\right) ,-\left[\bullet ,\bullet
\right] \right) \longleftrightarrow \left( \mathcal{F}\left( \mathcal{M}%
\right) ,\left\{ \bullet,\bullet \right\} ^c\right)  \label{iso1}
\end{equation}%
between the Lie algebras of real smooth functions and contact vector fields. 

According to \eqref{L-X-eta}, the flow of a contact Hamiltonian system preserves the contact structure, but it does not preserve neither the contact one-form nor the Hamiltonian function. Instead we obtain
\begin{equation}
{\mathcal{L}}_{X^c_H} \, H = - \mathcal{R}(H) H.
\end{equation}
Being a non-vanishing top-form we can consider $d\eta^n
\wedge \eta$ as a volume form on $\mathcal{M}$.  
Hamiltonian motion does not preserve the volume form since
\begin{equation}
{\mathcal{L}}_{X^c_H}  \, (d\eta^n
\wedge \eta) = - (n+1)  \mathcal{R}(H) d\eta^n
\wedge \eta.
\end{equation}%
However, it is immediate to see that, for a nowhere vanishing Hamiltonian function $H$, the quantity $ {H}^{-(n+1)}   (d\eta)^n \wedge\eta$ 
is preserved along the motion (see \cite{BrLeMaPa20}).

Referring to the Darboux' coordinates $(q^i,p_i,z)$, for a Hamiltonian function $H$, the Hamiltonian vector field, determined in \eqref{contact}, is computed to be
\begin{equation}\label{con-dyn}
X^c_H=\frac{\partial H}{\partial p_i}\frac{\partial}{\partial q^i}  - \big (\frac{\partial H}{\partial q^i} + \frac{\partial H}{\partial z} p_i \big)
\frac{\partial}{\partial p_i} + (p_i\frac{\partial H}{\partial p_i} - H)\frac{\partial}{\partial z},
\end{equation}
whereas the contact Jacobi bracket \eqref{cont-bracket} is 
\begin{equation}\label{Lag-Bra}
\{F,H\}^c = \frac{\partial F}{\partial q^i}\frac{\partial H}{\partial p_i} -
\frac{\partial F}{\partial p_i}\frac{\partial H}{\partial q^i} + \big(F  - p_i\frac{\partial F}{\partial p_i} \big)\frac{\partial H}{\partial z} -
\big(H  - p_i\frac{\partial H}{\partial p_i} \big)\frac{\partial F}{\partial z}.
\end{equation}
So, we obtain that the Hamilton's equations for $H$ as
\begin{equation}\label{conham}
\dot{q}^i= \frac{\partial H}{\partial p_i}, \qquad \dot{p}_i = -\frac{\partial H}{\partial q^i}- 
p_i\frac{\partial H}{\partial z}, \quad \dot{z} = p_i\frac{\partial H}{\partial p_i} - H.
\end{equation}

\textbf{Evolution vector fields}
Another vector field can be defined from a Hamiltonian function $H$ on a contact manifold $(M,\mathcal{\eta})$: the \emph{evolution vector field} of $H$~\cite{simoes2020contact}, denoted as $\varepsilon_H$, which is the one that satisfies
\begin{equation}\label{evo-def} 
\mathcal{L}_{\varepsilon_H}\eta=dH-\mathcal{R}(H)\eta,\qquad \eta(\varepsilon_H)=0.
\end{equation} 
In local coordinates it is given by
\begin{equation}\label{evo-dyn}
	\varepsilon_H=\frac{\partial H}{\partial p_i}\frac{\partial}{\partial q^i}  - \big (\frac{\partial H}{\partial q^i} + \frac{\partial H}{\partial z} p_i \big)
	\frac{\partial}{\partial p_i} + p_i\frac{\partial H}{\partial p_i} \frac{\partial}{\partial z},
\end{equation}
so that the integral curves satisfy the evolution equations
\begin{equation}\label{evo-eq}
	\dot{q}^i= \frac{\partial H}{\partial p_i}, \qquad \dot{p}_i = -\frac{\partial H}{\partial q^i}- 
	p_i\frac{\partial H}{\partial z}, \quad \dot{z} = p_i\frac{\partial H}{\partial p_i}.
\end{equation}

The evolution and Hamiltonian vector fields are related by
\begin{equation}
	\varepsilon_H = X^c_H + H \mathcal{R}.
\end{equation}

 \textbf{Quantomorphisms.} By asking the conformal factor $\mu$ in the definition \eqref{Cont-Dif} to be the unity, one arrives the conservation of the contact forms 
  \begin{equation}
        \varphi^*\eta_2 = \eta_1.
    \end{equation}
We call such a mapping as a strict contact diffeomorphism (or quantomorphism). For a contact manifold $(\mathcal{M},\eta)$
we denote the space of all strict contact transformations  as
\begin{equation}
{\rm Diff}_{con}^{st} ( \mathcal{M})
 =\left\{ \varphi \in {\rm Diff} ( \mathcal{%
M} ) :\varphi ^{\ast }\eta =   \eta \right\} \subset {\rm Diff}_{con} ( \mathcal{M})
.
\end{equation}
The Lie algebra of this group is consisting of the infinitesimal quantomorphisms
\begin{equation}\label{X-st-con}
\mathfrak{X}_{con}^{st} ( \mathcal{M} ) =\left\{ X\in \mathfrak{%
X}_{con} ( \mathcal{M} ) :\mathcal{L}_{X_{H}}\eta =0\right\}. 
\end{equation}
If the contact vector field is determined through a smooth function $H$ as in \eqref{contact} then $X_H$ falls into the subspace $\mathfrak{X}_{con}^{st} ( \mathcal{M} )$ if and only if $\lambda=-dH(\mathcal{R})=0$. This reads that, to generate an infinitesimal quantomorphism a function $H$ must not depend on the fiber variable $z$. 

Now consider the canonical contact manifold $(\mathcal{T}^*\mathcal{Q}, \eta_\mathcal{Q})$. For two functions those are not dependent on the fiber variable $z$, the contact Jacobi bracket $\{\bullet ,\bullet \}^c$ in \eqref{Lag-Bra} locally turns out to be equal to the canonical Poisson bracket on $T^*Q$, therefore we have that 
\begin{equation}
\left[ X_{H}^{st},X_{F}^{st}\right]=-X^{st}_{\left\{ H,F\right\} },
\end{equation}
where $ X_{H}^{st}$ is the infinitesimal quantomorphism generated by $H$. 
Accordingly, one arrives at an isomorphism
\begin{equation}
 \mathfrak{X}_{ham} (
T^{\ast }\mathcal{Q} ) \longrightarrow  \mathfrak{X}_{con}^{st} ( \mathcal{T}^*\mathcal{Q} ) , \qquad 
X_H\mapsto X_H-H\frac{\partial}{\partial z}
\end{equation}
from the Lie algebra of Hamiltonian vector
fields $\mathfrak{X}_{ham} (
T^{\ast }\mathcal{Q} )$ into the canonical cotangent bundle $T^{\ast }\mathcal{Q}$ to  the Lie algebra of infinitesimal quantomorphisms $ \mathfrak{X}_{con}^{st} ( \mathcal{T}^*\mathcal{Q})$ on the extended cotangent bundle.

\subsection{Contact Lagrangian Dynamics} 

Once more, we consider the extended configuration space $\mathcal{Q}\times \mathbb{R}$ but this time take it as the total space of the standard fibration from $\mathcal{Q}\times \mathbb{R}$ to $\mathbb{R}$. In this bundle structure, the base manifold is $\mathbb{R}$ and the fibration is simply the projection to the second factor. The first jet manifold is diffeomorphic to
\begin{equation}
\mathcal
{T}\mathcal{Q}=T\mathcal{Q}\times \mathbb{R}.
\end{equation}
We call this space as the extended tangent bundle. 

Suppose that $L:\mathcal
{T}\mathcal{Q}\mapsto \mathbb{R}$ is a Lagrangian function. In order to arrive at the dynamical equations governed by such a Lagrangian function, one needs to employ Herglotz principle  which is defined by an action functional \cite{Guenther,Herglotz,de2020review}. 

The value of the
functional attains its extremum if  ${q}(t)$ is a solution of the \emph{Herglotz equations} (also known as the generalized
Euler-Lagrange equations):
\begin{equation}\label{Herglotz}
\dot{q}^i=\frac{dq^i}{dt},\qquad 
\frac{\partial L}{\partial q^i} - \frac{d}{dt}\Big(\frac{\partial L}{\partial {\dot q}^i} \Big)
+ \frac{\partial L}{\partial z}\frac{\partial L}{\partial {\dot q}^i} = 0, 
\end{equation}
and $z$ is a solution of the Cauchy problem
\begin{equation}
\label{herglotzprinciple} \dot{z} = L(t,q^i,\dot{q}^i,z), \qquad 0 \leq t \leq \tau.
\end{equation}
It is important to notice that (\ref{Herglotz}) represents a family of 
differential equations since for each
function $q(t)$ a different differential equation arises, hence $z(t)$ depends on $q(t)$.
Without the explicit dependence of $z$, this problem reduces to a problem of the classical calculus of variations. If the functional $z$ defined in \eqref{herglotzprinciple} is invariant 
with respect to translation in time, then the quantity
\begin{equation} \label{consrule}
I = exp\Big(- \int^t  \frac{\partial L}{\partial z} d\theta \Big)\Big(L ( q, \dot{q}, z) - 
\frac{\partial L}{\partial \dot{q}^i}\dot{q}^i \Big)
\end{equation}
is conserved on solutions of the Herglotz equations for regular Lagrangians. 

For a regular Lagrangian function $L$, the fiber derivative determines a diffeomorphism from the extended tangent bundle $\mathcal
{T}\mathcal{Q}$ to the extended cotangent bundle $\mathcal
{T}^*\mathcal{Q}$ as
\begin{equation}\label{Leg-Trf}
\mathbb{F}L^c: \mathcal
{T}\mathcal{Q} \longrightarrow  \mathcal
{T}^*\mathcal{Q}, \qquad (q^i,\dot{q}^i,z)\mapsto 
(q^i,\frac{\partial L}{\partial \dot{q}^j},z)
\end{equation}
A direct calculation shows that the Legendre transformation~\eqref{Leg-Trf} maps the Herglotz equations in~\eqref{Herglotz} to the contact Hamilton's equations \eqref{conham} if the Hamiltonian function is defined to be
\begin{equation}
H(q^i,p_i,z)=\dot{q}^ip_i-L(q,\dot{q},z).
\end{equation} 

\textbf{The evolution-Herglotz equations}.
One can obtain the Lagrangian formalism for the evolution vector field by using a nonlinear nonholonomic action principle, as shown in~\cite{simoes2020geometry}. The resulting equations are the evolution Herglotz equations
\begin{equation}\label{Herglotz-evo-}
	\begin{gathered}
	\frac{\partial L}{\partial q^i} - \frac{d}{dt}\Big(\frac{\partial L}{\partial {\dot q}^i} \Big)
	+ \frac{\partial L}{\partial z}\frac{\partial L}{\partial {\dot q}^i} = 0, \\
	\dot{z} = \dot{q}^i \frac{\partial L}{\partial \dot{q}^i}.
	\end{gathered}
\end{equation}
	The Legendre transformation~\eqref{Leg-Trf} maps the evolution Herglotz equations in \eqref{Herglotz-evo-} to the evolution contact Hamiltonian dynamics in \eqref{evo-eq}.
 
\section{Tulczyjew's triple for Contact Geometry}\label{Section-TTc}

\subsection{Special Contact Structures} 

	Let $(\mathcal{M},\eta)$ be a contact manifold and the total space of a fibre bundle $(\mathcal{M},\rho,\mathcal{Q})$. We introduce a special contact structure as a quintuple
		\begin{equation} \label{scs-5}  
		(\mathcal{M},\rho,\mathcal{Q},\eta,\Phi),
				\end{equation}
where $
		\Phi$ is a fiber preserving contact
		diffeomorphism from $\mathcal{M}$ to the canonical contact manifold $(\mathcal{T}^*\mathcal{Q},\eta_{\mathcal{Q}})$. The two-tuple $(\mathcal{M},\eta)$ is said to be the underlying contact manifold of the special contact structure. Here, we have a diagram exhibiting a   special contact structure in a pictorial way
		\begin{equation} \label{scs}
		\xymatrix{ \mathcal{T}^*\mathcal{Q}\ar[ddr]_{\pi^0_\mathcal{Q}} &&\mathcal{M}
			\ar[ddl]^{\rho} \ar[ll]_{\Phi}
			\\ \\&\mathcal{Q} }  
		\end{equation}
where $\pi^0_\mathcal{Q}$ is the fibration given in \eqref{Sec-Jet}. 

It is possible to define a special contact space starting with a special symplectic space defined in Section \ref{Ss-sss}. 
For this, suppose that $(\mathcal{P},\omega=-d\theta)$ is an exact symplectic manifold and in the product manifold $\mathcal{P}\times \mathbb{R}$ we consider the standard contact structure $\eta=dz-\theta$ (that is, the contactization of the exact symplectic structure $\omega=-d\theta$). Assume also that $\mathcal{P}$ admits a special symplectic structure $		(\mathcal{P},\pi,\mathcal{Q},\theta,\phi)$ as pictured in \eqref{sss-}. Then, $\mathcal{M}$ admits a special contact structure $  		(\mathcal{M},\rho,\mathcal{Q},\eta,\Phi)$, where $\Phi(p,z)=(\phi(p),z)$, and the following diagram
			\begin{equation} \label{scs-3}
		\xymatrix{\mathcal{T}^*\mathcal{Q}
			\ar[dd]_{\pi_\mathcal{Q}^1} &&\mathcal{M}=\mathcal{P}\times \mathbb{R} \ar[ll]_{\Phi}\ar[dd]^{\rm pr} 
\\ \\		T^*\mathcal{Q}
			\ar[ddr]_{\pi_{\mathcal{Q}}} 	
		&&\mathcal{P} \ar[ll]^{\phi}\ar[ddl]^{\pi}
		\\	\\&\mathcal{Q} }  
		\end{equation}	
is commutative, with $\pi^1_\mathcal{Q}$ the fibration given in \eqref{Sec-Jet}. 		Here, it is considered that $\rho=\pi \circ {\rm pr}$. This construction is the contactization of the special symplectic structure. 

We now merge a Morse family $E$ defined on a fiber bundle $(\mathcal{W},\tau,\mathcal{Q})$ and a special contact space $  		(\mathcal{M},\rho,\mathcal{Q},\eta,\Phi)$ in order to arrive at a Legendrian submanifold of $(\mathcal{M},\eta)$. For this, consider the following commutative diagram
\begin{equation} \label{Morse-Gen--}
\xymatrix{
	\mathbb{R}& \mathcal{W} \ar[dd]^{\tau}\ar[l]^{E}& \mathcal{T}^*\mathcal{Q} \ar@(ul,ur)^{ \mathcal{N}} \ar[ddr]_{\pi_{\mathcal{Q}}^0}& &\mathcal{M}\ar@(ul,ur)^{ \mathcal{N}_{E}} \ar[ll]_{\Phi} \ar[ddl]^{\rho}\\ \\ &
	\mathcal{Q} \ar@{=}[rr]&& \mathcal{Q} 
}
\end{equation}  
Referring to the definition in \eqref{MFGen-C}, we obtain a Legendrian submanifold $\mathcal{N}$ of the jet bundle $\mathcal{T}^*\mathcal{Q}$. Then, by employing the inverse of the contact diffeomorphism $\Phi$, we arrive at a Legendrian submanifold $\mathcal{N}_{E}$ of $\mathcal{M}$. Referring to this realization, we shall exhibit both the contact Hamiltonian and contact Lagrangian dynamics as Legendrian submanifolds of the same contact manifold in the following subsection.

\subsection{Tangent Contact Manifold} \label{Sec-Tcont}

We start by lifting a contact structure $\eta$ on a contact manifold $\mathcal{M}$ a contact structure on the extended tangent bundle $\mathcal{T}\mathcal{M}$. This lifting is in introduced in \cite{IbLeMaMa97} to characterize the contact vector fields on $\mathcal{M}$ (in particular, the Hamiltonian vector fields in $\mathcal{M}$) in terms of Legendrian submanifolds of the contact manifold $\mathcal{T}\mathcal{M}$. In fact, more later, in this direction, we shall use some other results those available in \cite{IbLeMaMa97}. 
\begin{theorem}
For a contact manifold $(\mathcal{M},\eta)$, the extended tangen bundle $\mathcal{T}\mathcal{M} \simeq T\mathcal{M} \times \mathbb{R}$ is a contact manifold by admitting a contact one-form
\begin{equation}\label{eta-T-}
\eta^\mathcal{T} := u {\eta}^V+{\eta}^C
\end{equation}
where $u$ is coordinate on $\mathbb{R}$ whereas ${\eta}^C$ and ${\eta}^V$ are the complete and vertical lifts of $\eta$, respectively. 
\end{theorem} 
The one-form $\eta^\mathcal{T}$ is said to be the tangent contact structure and we will denote the tangent contact manifold as a two-tuple
\begin{equation}\label{T-cont}
(\mathcal{T}\mathcal{M},\eta^\mathcal{T})=( T\mathcal{M} \times \mathbb{R}, u {\eta}^V+{\eta}^C).
\end{equation}

\textbf{Contact Hamiltonian Dynamics as a Legendrian Submanifold.}    Let $(\mathcal{M},\eta)$ be a contact manifold. Consider a vector field $X$, a real valued function $\lambda$ on $\mathcal{M}$, hence a section 
\begin{equation}\label{bar-x}
            (X,\lambda): \mathcal{M} \longrightarrow  \mathcal{T}\mathcal{M}=T\mathcal{M}\times \mathbb{R}, \qquad 
            m \mapsto (X(m), \lambda(m)).
\end{equation}
of the fibration $\tau^0_\mathcal{M}:\mathcal{T}\mathcal{M} \mapsto \mathcal{M}$. We plot the following commutative diagram to see this 
		\begin{equation} \label{scs-4}
		\xymatrix{ \mathcal{T}\mathcal{M}  \ar[rrd]^{\tau^1_\mathcal{M}} \ar[ddr]^{\tau^0_\mathcal{M}} 
			\\  && T\mathcal{M}
			\ar[dl]_{\tau_{\mathcal{M}}}  \\&\mathcal{M}
			\ar@/^1pc/[uul]^{(X,\lambda)} \ar@/_1pc/[ur]_{X}
			 }  
		\end{equation}
Using Theorem 3.13 in  \cite{IbLeMaMa97} and the comments at the beginning of this subsection, we deduce that the pair 
 $(X,\lambda)$ is a contact vector field (an infinitesimal conformal contactomorphism), that is, an element of $\mathfrak{X}_{con} ( \mathcal{M} ) $ in \eqref{algcon} if and only if the image space of $(X,\lambda)$ is a Legendrian submanifold of the tangent contact manifold $(\mathcal{T}\mathcal{M},\eta^\mathcal{T})$. This result states that the image of a contact Hamiltonian vector field $X^c_H$, after suitably included in the contactified tangent bundle, turns out to be a Legendrian submanifold. As discussed in the previous section, the conformal factor $\lambda$ in the present case is $\mathcal{R}(H)$. So that, the image of the mapping 
\begin{equation}\label{bar-x-H}
(X^c_H, \mathcal{R}(H)):\mathcal{M} \longrightarrow \mathcal{T}\mathcal{M}=T\mathcal{M}\times \mathbb{R}, \qquad m\mapsto (X^c_H(m),\mathcal{R}(H)(m))
\end{equation}
is a Legendrian submanifold of  the tangent contact manifold $\mathcal{T}\mathcal{M}$. 
To see this more clearly, let us discuss this geometry in the realm of special contact spaces. 

Consider a contact manifold $(\mathcal{M},\eta)$. Its extended 
tangent bundle $\mathcal{T}\mathcal{M}$ is a contact manifold
endowed with the contact structure given by~\eqref{T-cont}  
and $\mathcal{T}^*\mathcal{M}$, as an extended cotangent bundle, also is a contact manifold. Between these two extended spaces, we introduce a fiber preserving contact diffeomorphism $\beta^c$ as follows
\begin{equation}\label{beta-c}
\beta^c:\mathcal{T}\mathcal{M}  \longrightarrow \mathcal{T}^*\mathcal{M},\qquad (V,u) \mapsto \big (-\iota _ V d \eta - u \eta, \eta(V)\big)
\end{equation}
where $\eta$ is the contact one-form on $\mathcal{M} $.  It is a direct computation to see that
\begin{equation}\label{beta-c-pf}
(\beta^c)^*\eta_{\mathcal{T}^*\mathcal{M}}=\eta^\mathcal{T}
\end{equation}
where $\eta_{\mathcal{T}^*\mathcal{M}}$ is the contact one-form on the canonical contact manifold $\mathcal{T}^*\mathcal{M}$ and $\eta^\mathcal{T}$ is the lifted contact one-form on $\mathcal{T}\mathcal{M} $ given in \eqref{eta-T-}. This observation permits us to determine a special contact space and following the order given in \eqref{scs-5}, we write this special contact manifold as
\begin{equation}\label{scs-1}
(\mathcal{T}\mathcal{M},\tau^0_\mathcal{M},\mathcal{M},\eta^\mathcal{T},\beta^c).
\end{equation}
Accordingly, following the picture in \eqref{scs}, we plot the following diagram 		
\begin{equation} \label{t-scs}
		\xymatrix{\mathcal{T}\mathcal{M}
			\ar[ddr]^{\tau^0_{\mathcal{M}}} \ar[rr]^{\beta^c} &&\mathcal{T}^*\mathcal{M} \ar[ddl]_{\pi^0_{\mathcal{M}}}
			\\ \\&\mathcal{M} \ar@/^1pc/[uul]^{(X^c_H,\mathcal{R}(H))} \ar@/_1pc/[uur]_{-\mathcal{T}^* H}} 
		\end{equation}
		where we have employed the projections $\pi^0_{\mathcal{M}}:\mathcal{T}^*\mathcal{M}\mapsto \mathcal{M}$ and  $\tau^0_{\mathcal{M}}:\mathcal{T}\mathcal{M}\mapsto \mathcal{M}$. This diagram also manifesting how one can transfer the Legendrian submanifolds one onto the other. If $H$ is a Hamiltonian function on the contact manifold $\mathcal{M}$, then the image space of 
		$-\mathcal{T}^* H$, as defined in \eqref{j1F}, is a Legendrian submanifold of $\mathcal{T}^*\mathcal{M}$. Thus, since  $\beta^c$ is a contact diffeomorphism, by pulling back the image space of $-\mathcal{T}^* H$, we arrive at a Legendrian submanifold of the tangent contact manifold $\mathcal{T}\mathcal{M}$. But, using \eqref{beta-c}, we have that
\begin{equation}\label{contact-Ham-Leg}
\beta^c\circ (X^c_H,\mathcal{R}(H)) = -\mathcal{T}^* H=-(dH,H).
\end{equation}
This is the contact version of the identity \eqref{bemol-H}. Notice that, this observation can be considered as an indirect proof of the assertion that the image space of a contact vector field is a Legendrian submanifold.   

\textbf{Local Picture.}
Consider Darboux' coordinates $(q^i,p_i,z)$ on $ \mathcal{M} $ then we assume the induced coordinates on $\mathcal{T}\mathcal{M}\simeq T\mathcal{M} \times \mathbb{R}$ as $(q^i,p_i,z,\dot{q}^i,\dot{p}_i,\dot{z},u)$. In this local realization, the lifted contact one-form $ \eta^\mathcal{T}$ defined in \eqref{eta-T-} is computed to be
\begin{equation}\label{eta-T}
\begin{split}
 \eta^\mathcal{T} &=  u{\eta}^V + \eta^C 
  =          u (\dd z - p_i \dd q^i) + (\dd \dot{z} - \dot{p}_i\dd q^i - p_i \dd \dot{q}^i)
  \\ 
  &=\dd \dot{z} +u dz -(\dot{p}_i+u p_i) \dd q^i -  p_i \dd \dot{q}^i,
  \end{split}
        \end{equation}
and the Reeb vector field is $\mathcal{R}^{\mathcal{T}}=\partial/\partial \dot{z}$. In this realization, the contact mapping $\beta^c$ in \eqref{beta-c} turns out to be
\begin{equation}\label{beta-c-}
\beta^c(q^i,p_i,z,\dot{q}^i,\dot{p}_i,\dot{z}, u) =(q^i,p_i,z, u p_i+\dot{p}_i ,- \dot{q}^i, -u,\dot{z}-p_i \dot{q}^i).
\end{equation} 
The fact that $\beta^c$ is a contactomorphism, that is the identity \eqref{beta-c-pf}, follows from a direct calculation in coordinates. Observe that,
\begin{equation}
\begin{split}
        (\beta^c)^*(\eta_{\mathcal{T}^* \mathcal{M}}) &= 
        \dd (\dot{z}-p_i \dot{q}^i)-(u p_i+\dot{p}_i)q^i+\dot{q}^i \dd p_i +u\dd z
         \\ &= u (\dd z - p_i \dd q^i) + (\dd \dot{z} - \dot{p}_i\dd q^i - p_i \dd \dot{q}^i)=
       \eta^\mathcal{T}.
\end{split}
 \end{equation}
Consider now a Hamiltonian function $H$ on the contact manifold $\mathcal{M}$. Minus of its first prolongation defines a Legendrian submanifold of $\mathcal{T}^*\mathcal{M}$ given by
\begin{equation}\label{Ham-Leg-Sub}
{\rm im}(-\mathcal{T}^*H)=\big\{(q^i,p_i,z, -\frac{\partial H}{\partial q^i},-\frac{\partial H}{\partial p_i},-\frac{\partial H}{\partial z},-H )\in \mathcal{T}^*\mathcal{M}:H=H(q,p,z) \big\}. 
\end{equation}
 Referring to the inverse of the contact diffeomorphism $\beta^c$ we compute a Legendrian submanifold of $\mathcal{T}\mathcal{M}$ as
\begin{equation}\label{N_H} 
\begin{split}
\mathcal{N}_{-H}&= (\beta^c)^{-1}\big({\rm im}(-\mathcal{T}^*H)\big)\\
&= \big\{(q^i,p_i,z, \frac{\partial H}{\partial p_i}, -\frac{\partial H}{\partial z}p_i-\frac{\partial H}{\partial q^i},
p_i\frac{\partial H}{\partial p_i}-H,\frac{\partial H}{\partial z})\in \mathcal{T}\mathcal{M}:H=H(q,p,z)   \big\}. 
\end{split}
 \end{equation}
 It is immediate to see that the Legendrian submanifold can be alternatively  obtained by 
\begin{equation} 
 \mathcal{N}_{-H}=\im (X^c_H,\mathcal{R}(H)).
  \end{equation} 
So, the contact Hamiltonian dynamics is determined by the Legendrian submanifold $\mathcal{N}_{-H}$ as follows. Let $\sigma:I\mapsto\mathcal{M}$ be a smooth curve on $\mathcal{M}$ and consider the lift of $\sigma$ to $\mathcal{T}\mathcal{M}$  defined by 
\begin{equation}\label{sigma-T}
 \sigma^{\mathcal{T}}=(\dot{\sigma}(t), \mathcal{R}(H)(\sigma(t)),\qquad \forall t\in I,
  \end{equation}
where $\dot{\sigma}$ is the tangent lift to  $T\mathcal{M}$.
Then $\sigma$ is a solution of the contact Hamilton's equations for $H$ if and only if its lift $\sigma^{\mathcal{T}}$ to $\mathcal{T}\mathcal{M}$ is contained in the Legendrian submanifold $\mathcal{N}_{-H}$. In fact, if the local expression of $\sigma$ is
\begin{equation}
  \sigma(t)=(q^i(t),p_i(t),z(t))
    \end{equation}
then, from \eqref{N_H}, it follows that  $\sigma^{\mathcal{T}}\in \mathcal{N}_{-H}$ for every $t\in I$, if and only if $\sigma$ satisfies the following equations
\begin{equation}\label{conHamee}
\frac{dq^i}{dt}=\frac{\partial H}{\partial p_i},\qquad \frac{dp_i}{dt}=-\frac{\partial H}{\partial z}p_i-\frac{\partial H}{\partial q^i}, \qquad \frac{dz}{dt} =p_i\frac{\partial H}{\partial p_i}-H.
\end{equation}
They are precisely the contact Hamilton's equations in \eqref{conham}.

\subsection{Contact Tulczyjew's Triple} 

Classical Tulczyjew's triple is obtained by properly merging two special symplectic structures. Following the same understanding, we introduce Tulczyjew's triple for contact dynamics by properly merging two special contact structures.

\textbf{Contact Lift of Diffeomorphism $\psi$ in \eqref{shuttle-loc}.} Recall the symplectic diffeomorphism $\psi$, defined in \eqref{shuttle-loc}, from the cotangent bundle $T^*T\mathcal{Q}$ to the cotangent bundle $T^*T^*\mathcal{Q}$. We extend this mapping to the level of contact manifolds 
\begin{equation}
\mathcal{T}^*\mathcal{T}^*\mathcal{Q}\simeq T^*(T^*\mathcal{Q}\times\mathbb{R})\times\mathbb{R}, \qquad \mathcal{T}^*\mathcal{T}\mathcal{Q}
\simeq 
T^*(T\mathcal{Q}\times\mathbb{R})\times\mathbb{R}
\end{equation}
by assuming that extension is the identity on $T\RR$ and it vanishes on the zero section. Accordingly, we recall the mapping $\psi$ in \eqref{shuttle-loc}, consider a local chart $(q^i,\dot{q}^i,z,a_i,\dot{a}_i,v,u)$ on $\mathcal{T}^*\mathcal{T}\mathcal{Q} $ then we have an arbitrary extension
\begin{equation*} 
\psi^c:\mathcal{T}^*\mathcal{T}\mathcal{Q}  \longrightarrow \mathcal{T}^*\mathcal{T}^*\mathcal{Q}, \qquad 
(q^i,\dot{q}^i,z,a_i,\dot{a}_i,v,u)  \mapsto
        (q^i,\dot{a}_i,z,a_i,-\dot{q}^i,v,w).
\end{equation*}  
To determine $w$, we compute
\begin{equation}
  (\psi^c)^* \eta_{\mathcal{T}^*\mathcal{Q}} = d w - \psi^*\theta_{T^*\mathcal{Q}} - v d z = d w - d u  + d (\dot{q}^i \dot{a}_i) + \eta_{\mathcal{T}\mathcal{Q}},
\end{equation}
where $\eta_{\mathcal{T}^*\mathcal{Q}}$ and $\eta_{\mathcal{T}\mathcal{Q}}$ are the canonical contact one-forms on $\mathcal{T}^*\mathcal{T}^*\mathcal{Q}$ and $\mathcal{T}^*\mathcal{T}\mathcal{Q}$, respectively. 
Hence $\psi^c$ is a contact mapping if and only if $d w =  d(u - \dot{q}^i \dot{a}_i) $. Since, we ask $w$ vanish on the zero section, necessarily $ w = u - \dot{q}^i \dot{a}_i$. Therefore, the unique extension is, locally, computed to be
    \begin{equation}\label{psi-c}
        \psi^c: \mathcal{T}^*\mathcal{T}\mathcal{Q}  \longrightarrow \mathcal{T}^*\mathcal{T}^*\mathcal{Q}, \qquad 
       (q^i,\dot{q}^i,z,a_i,\dot{a}_i,v,u) \mapsto
        (q^i,\dot{a}_i,z,a_i,-\dot{q}^i,v,u- \dot{a}_i \dot{q}^i).
    \end{equation}
Next, we will present an intrinsic definition of $\psi^c$. For this purpose, we will use the following identifications 
    \begin{equation}
\mathcal{T}^*\mathcal{T}\mathcal{Q}\cong (T^*T\mathcal{Q}\times T^*\mathbb{R})\times \mathbb{R}, \qquad \mathcal{T}^*\mathcal{T}^*\mathcal{Q}\cong (T^*T^*\mathcal{Q}\times T^*\mathbb{R})\times \mathbb{R}.
    \end{equation}
    Then, we have that 
   \begin{equation}
   \psi^c((\gamma,vdz),u)=\big(\psi(\gamma),vdz,u-\langle (\pi_{T^*\mathcal{Q}}(\psi (\gamma)),\pi_{T\mathcal{Q}}(\gamma) \rangle \big)
    \end{equation}
    where $\langle \bullet, \bullet \rangle$ is the canonical pairing between $T^*\mathcal{Q}$ and $T\mathcal{Q}$.
    
\textbf{The Left Wing of the triple.}
Now, using the contact mapping $\psi^c$ in \eqref{psi-c}, we define a contact diffeomorphism
    \begin{equation}  \label{alpha-c}
    \begin{split}
    \alpha^c: \mathcal{T}\mathcal{T}^*\mathcal{Q}  \longrightarrow \mathcal{T}^*\mathcal{T}\mathcal{Q}, \qquad V&\mapsto (\psi^c)^{-1} \comp \beta^c(V)
    \\
    (q^i,p_i,z,\dot{q}^i,\dot{p}_i,\dot{z}, u) &\mapsto
          (q^i,\dot{q}^i,z,up_i+ \dot{p}_i, p_i,-u,\dot{z}).
    \end{split}
    \end{equation}
Indeed, it is immediate to show that 
 \begin{equation}
     (\alpha^c)^*\eta_{\mathcal{T}\mathcal{Q}}=\eta^{\mathcal{T}}_\mathcal{Q},
       \end{equation}
    where $\eta^{\mathcal{T}}_\mathcal{Q}$ is the lifted contact one-form on $\mathcal{T}\mathcal{T}^*\mathcal{Q}$ given in \eqref{eta-T}.
Then, we define the following special contact structure 
\begin{equation}\label{scs-2-set}
(\mathcal{T}\mathcal{T}^*\mathcal{Q},\mathcal{T}\pi^0_{\mathcal{Q}},\mathcal{T}\mathcal{Q}, \eta_{\mathcal{Q}}^\mathcal{T},\alpha^c)
\end{equation}
which is diagrammatically given by
\begin{equation} \label{scs-2}
\xymatrix{\mathcal{T}^*\mathcal{T}\mathcal{Q} \ar[ddr]^{\pi^0_{\mathcal{T}\mathcal{Q}}} && \mathcal{T}\mathcal{T}^*\mathcal{Q}
\ar[ddl]^{\mathcal{T}\pi^0_{\mathcal{Q}}} \ar[ll]_{\alpha^c}
	\\\\&\mathcal{T}\mathcal{Q}\ar@/^1pc/[uul]^{\mathcal{T}^* L}}  
\end{equation}
with the projections
\begin{equation}\label{proj-set-2}
\begin{split}
\pi^0_{\mathcal{T}\mathcal{Q}}&:\mathcal{T}^*\mathcal{T}\mathcal{Q}\simeq T^*\mathcal{T}\mathcal{Q}\times \mathbb
{R} 
\longrightarrow  \mathcal{T}\mathcal{Q},\qquad (W,u)\mapsto \pi_{\mathcal{T}\mathcal{Q}}(W)
\\
\mathcal{T}\pi^0_{\mathcal{Q}}&:\mathcal{T}\mathcal{T}^*\mathcal{Q} \simeq T \mathcal{T}^*\mathcal{Q}\times \mathbb
{R} 
\longrightarrow  \mathcal{T}\mathcal{Q},\qquad (V,s)\mapsto (T\pi_{\mathcal{Q}}(U),z).
\end{split}
 \end{equation}
Here, we have employed the following global trivialization $V=(U,z,\dot{z})$ in $TT^*\mathcal{Q}\times T\mathbb{R}$.
 
\textbf{Generalized Euler-Lagrange Equations as  a Lagrangian Submanifold.}  Consider a Lagrangian function $L$ on $\mathcal{T}\mathcal{Q}=T\mathcal{Q}\times \mathbb{R}$. The image space of its first prolongation, that is, ${\rm im}(\mathcal{T}^*L)$, is a Legendrian submanifold of $\mathcal{T}^*\mathcal{T}\mathcal{Q}$ computed to be
\begin{equation}
{\rm im}(\mathcal{T}^*L)=\big\{(q^i,\dot{q}^i,z, \frac{\partial L}{\partial q^i},\frac{\partial L}{\partial \dot{q}^i},\frac{\partial L}{\partial z},L)\in \mathcal{T}^*\mathcal{T}\mathcal{Q} : L=L(q,\dot{q},z)\big\}\subset \mathcal{T}^*\mathcal{T}\mathcal{Q}
\end{equation}
Referring to the left wing of the contact triple \eqref{con-TT}, that is by applying the inverse of the mapping $\alpha^c$, we  arrive at a Legendrian submanifold of $\mathcal{T}\mathcal{T}^*\mathcal{Q}$ as
\begin{equation}\label{N_L}
\begin{split}
\mathcal{N}_L&= (\alpha^c)^{-1}({\rm im}(\mathcal{T}^*L)) \\
&=\big\{(q^i,\frac{\partial L}{\partial \dot{q}^i},z, \dot{q}^i, \frac{\partial L}{\partial z}\frac{\partial L}{\partial \dot{q}^i}+\frac{\partial L}{\partial q^i},L, -\frac{\partial L}{\partial z})\in \mathcal{T}\mathcal{T}^*\mathcal{Q} : L=L(q,\dot{q},z)\big\}\subset \mathcal{T}\mathcal{T}^*\mathcal{Q}.
\end{split}
\end{equation}
The Lagrangian dynamics is determined by the Legendrian submanifold $\mathcal{N}_L$ as follows. First of all, we consider the transformation $\mathbb{F}L:\mathcal{T}\mathcal{Q}\mapsto \mathcal{T}^*\mathcal{Q}$ induced by $L$ given by
\begin{equation}
\mathbb{F}L(u,z)=(\mathbb{F}L_1(u,z),z),
\end{equation}
with $\mathbb{F}L_1(u,z)$ defined by
\begin{equation}
\langle \mathbb{F}L_1(u,z),u' \rangle=\frac{d}{dt}\Big\vert_{t=0} L(u+tu',z).
\end{equation}
Then, a smooth curve $c:I\mapsto \mathcal{T}\mathcal{Q}$ is a solution of the dynamics if and only if the lift $(\mathbb{F}L\circ c)^{\mathcal{T}}:I\mapsto \mathcal{T}\mathcal{T}^*\mathcal{Q}$ of the curve $\mathbb{F}L\circ c:I \mapsto \mathcal{T}^*\mathcal{Q}$ given by
\begin{equation}
(\mathbb{F}L\circ c)^{\mathcal{T}}(t)=\big(\dot{\mathbb{F}L\circ c}(t),\frac{\partial L}{\partial z}(c(t))\big), \forall t\in I 
\end{equation}
is contained in the Legendrian submanifold  $\mathcal{N}_L$. In fact, if the local expression of $c$ is 
\begin{equation}
c(t)=(q^i(t),\dot{q}^i(t),z(t))
\end{equation}
then, using \eqref{N_L}, it follows that $(\mathbb{F}L\circ c)^{\mathcal{T}}(t)\in \mathcal{N}_L$, $ \forall t\in I$, if and only if $c$ satisfies 
\begin{equation}
\frac{dq^i}{dt}=\dot{q}^i,\qquad \frac{d}{dt}\big( \frac{\partial L}{\partial \dot{q}^i}\big)=\frac{\partial L}{\partial z}\frac{\partial L}{\partial \dot{q}^i}+\frac{\partial L}{\partial q^i}, \qquad \frac{dz}{dt}=L(q,\dot{q},z) 
\end{equation}
where the first two equations are the Herglotz equations in \eqref{Herglotz} and the latter is the Herglotz's differential principle in \eqref{herglotzprinciple}. We remark that  this theory does not require any regularity conditions on the Lagrangian function. 
 
 \textbf{The Right Wing of the triple.} To have the right wing, we simply replace the arbitrary contact manifold $(\mathcal{M},\eta)$  with the canonical contact manifold $(\mathcal{T}^*\mathcal{Q},\eta_\mathcal{Q})$ in the special contact manifold  \eqref{scs-1}.
In this case, the total spaces become $\mathcal{T}^*\mathcal{T}^*\mathcal{Q}$ and $\mathcal{T}\mathcal{T}^*\mathcal{Q}$ with the fibrations
\begin{equation}\label{proj-set-1}
\begin{split}
\tau^0_{\mathcal{T}^*\mathcal{Q}}&:\mathcal{T}\mathcal{T}^*\mathcal{Q}\simeq T\mathcal{T}^*\mathcal{Q}\times \mathbb
{R} 
\longrightarrow  \mathcal{T}^* \mathcal{Q},\qquad (V,u)\mapsto \tau_{\mathcal{T}^* \mathcal{Q}}(V),
\\
\pi^0_{\mathcal{T}^*\mathcal{Q}}&:\mathcal{T}^*\mathcal{T}^*\mathcal{Q} \simeq T^*\mathcal{T}^*\mathcal{Q}\times \mathbb
{R} 
\longrightarrow  \mathcal{T}^* \mathcal{Q},\qquad (Z,u)\mapsto \pi_{\mathcal{T}^* \mathcal{Q}}(Z).
\end{split}
 \end{equation}
We denote this special contact manifold as 
\begin{equation}\label{scs-1-set}
(\mathcal{T}\mathcal{T}^*\mathcal{Q},\tau^0_{\mathcal{T}^*\mathcal{Q}},\mathcal{T}^*\mathcal{Q},\eta_{\mathcal{Q}}^\mathcal{T},\beta^c).
\end{equation} 
Since Darboux' coordinates are employed in Subsection \ref{Sec-Tcont}, all the local formulations available in that subsection hold also in the present case.

\textbf{Contact  triple.} We now merge the special contact spaces \eqref{scs-2-set} and \eqref{scs-1-set} in order to construct a Tulzyjew's triple for the contact geometry. Accordingly, we couple the commutative diagrams in \eqref{t-scs} and \eqref{scs-2}
and arrive at the following diagram
\begin{equation} \label{con-TT}
\xymatrix{ \mathcal{T}^*\mathcal{T}\mathcal{Q}\ar[ddr]_{\pi^0_{\mathcal{T}\mathcal{Q}}} && \mathcal{T}\mathcal{T}^*\mathcal{Q}
\ar[ddl]^{\mathcal{T}\pi^0_{\mathcal{Q}}} \ar[ll]_{\alpha^c}\ar[rr]^{\beta^c}\ar[ddr]_{\tau^0_{\mathcal{T}^*\mathcal{Q}}}
&& \mathcal{T}^*\mathcal{T}^*\mathcal{Q}\ar[ddl]^{\pi^0_{\mathcal{T}^*\mathcal{Q}}} 
	\\\\&\mathcal{T}\mathcal{Q}\ar[ddr]_{\tau^0_{\mathcal{Q}}}  &&\mathcal{T}^*\mathcal{Q}\ar[ddl]^{\pi^0_{\mathcal{Q}}}
	\\\\
	&&\mathcal{Q}}
\end{equation}
where the contact diffeomorphisms $\beta^c$ and $\alpha^c$ are those defined in \eqref{beta-c} and \eqref{alpha-c}, respectively. Here, the projections $\tau^0_{\mathcal{T}^*\mathcal{Q}}$ and $\pi^0_{\mathcal{T}^*\mathcal{Q}}$ 
(respectively $\pi^0_{\mathcal{T}\mathcal{Q}}$ and $\mathcal{T}\pi^0_{\mathcal{Q}}$) are given by \eqref{proj-set-1} (respectively, \eqref{proj-set-2}).

\subsection{Evolution Contact Tulczyjew's Triple}

We will analyse how the Tulczyjew's triple can be used to understand the evolution dynamics. First of all, we notice that there is not, a natural way, to describe an evolution  vector field on a contact manifold $\mathcal{M}$ as a Legendrian submanifold of $\mathcal{T}\mathcal{M}$. However, we will see that such a vector field may be considered as a Lagrangian submanifold of an exact symplectic submanifold of $\mathcal{T}\mathcal{M}$.  

Consider a contact manifold $(\mathcal{M},\eta)$. In \eqref{decomp-TM}, we have realized the kernel of the contact one-form as a symplectic vector bundle $H(\mathcal{M})$ of the tangent bundle $T\mathcal{M}$. Consider now the inclusion mapping
	\begin{equation}\label{inc}
	j = (i,{\mathrm{Id}}_{\mathbb{R}}):H \mathcal{M} \times \mathbb{R} \hookrightarrow \mathcal{T}\mathcal{M} = T\mathcal{M} \times \mathbb{R}
	\end{equation}
	where $i$ is the inclusion of $H \mathcal{M}$ into $T\mathcal{M}$. 
Consider now the one-form $\theta_{\eta} = j^* \eta^\mathcal{T}$ on $H \mathcal{M}\times \mathbb{R}$, and the two-form $\omega_{\eta} =   d \theta_{\eta}$. In the following proposition, we state that $\omega_{\eta}$ is a symplectic two-form on $H(\mathcal{M})\times \mathbb{R}$
\begin{proposition}
The exact two-form $\omega_{\eta} =   d \theta_{\eta}$ induces a symplectic structure on $ H \mathcal{M}  \times \mathbb{R} $. 
\end{proposition}
{\rm\textbf{Proof.}} Let $(q^i,p_i,z)$ be the Darboux' coordinates on $\mathcal{M}$, so that  $(q^i,p_i,z, \dot{q}^i,\dot{p}_i,\dot{z}, u)$ are the induced  coordinates for $\mathcal{T}\mathcal{M}$. In the light of the inclusion \eqref{inc}, we can employ $(q^i,p_i,z, \dot{q}^i,\dot{p}_i, u)$ as a local coordinate chart on $H \mathcal{M}  \times \mathbb{R}$ since
	\begin{equation}\label{map-j}
		j(q^i,p_i,z, \dot{q}^i,\dot{p}_i, u) = (q^i,p_i,z, \dot{q}^i,\dot{p}_i, p_i \dot{q}^i, u).
	\end{equation}
In other words, 
\begin{equation}\label{HMxR}
H \mathcal{M}\times \mathbb{R}=\big \{ 
(q^i,p_i,z, \dot{q}^i,\dot{p}_i,\dot{z}, u)\in \mathcal{T}\mathcal{M}\times \mathbb{R} : \dot{z}-p_i\dot{q}^i=0
\big \}.
\end{equation}
So, the Reeb vector field $\mathcal{R}^{\mathcal{T}}=\partial/\partial \dot{z}$ of the contact manifold $(\mathcal{T}\mathcal{M},\eta^\mathcal{T})$ is transverse to the manifold of codimension one $H \mathcal{M}\times \mathbb{R}$. This implies that the exact two-form $\omega_{\eta} = d \theta_{\eta}=d(j^* \eta^\mathcal{T})$ induces a symplectic structure on  $H \mathcal{M}\times \mathbb{R}$. 

Note that, from \eqref{eta-T} and \eqref{map-j}, it follows that 
	\begin{equation}\label{omega-H}
	\begin{split}
		\theta_{\eta} &= 
		u dz -(\dot{p}_i+u p_i) \dd q^i + \dot{q}^i \dd  p_i,
\\
		\omega_{\eta}& = du \wedge d z - d \dot{p}_i \wedge \dd q^i - p_i du \wedge \dd q^i - u d p_i \wedge \dd q^i  + d \dot{q}^i \wedge d  p_i.
\end{split}
	\end{equation}

	\textbf{Lagrangian submanifolds and evolution vector fields.} Given a vector field $X$ on $\mathcal{M}$ and a real smooth function $f:\mathcal{M}\to \mathbb{R}$, one can construct a section $(X,f): \mathcal{M} \to \mathcal{T} \mathcal{M}$ of the extended tangent bundle $\mathcal{T} \mathcal{M}$.
	\begin{theorem}\label{lag-evo-char}
	Let $(\mathcal{M},\eta)$ be a contact manifold. 
	The map	$(X,f):\mathcal{M}\mapsto \mathcal{T}\mathcal{M}$ defines a Lagrangian submanifold of the exact symplectic manifold $(H \mathcal{M} \times \mathbb{R}, \omega_{\eta})$ if and only if $\eta(X) = 0$ and $\mathcal{L}_X \eta + f \eta$ is closed.
	\end{theorem}

	{\rm\textbf{Proof.}   }
	First of all, note that the image of $X$ lies on $H\mathcal{M}\times \mathbb{R}$ if and only if $\eta(X) = 0$.
	For the second condition we use well-known properties of complete and vertical lifts (see \cite{YaIs73});
	\begin{equation}
		(X, f)^* \eta^\mathcal{T} = 
		(X, f)^* (u \eta^V + \eta^C) = 
		f X^* \eta^V + X^* \eta^C = f \eta + \mathcal{L}_X \eta.
	\end{equation}
	Thus, the image of $(X, f)$ is Lagrangian if and only if
	\begin{equation}
		(X, f)^*\omega_{\eta} = d (X, f)^* \eta^\mathcal{T} =
		d( f \eta + \mathcal{L}_X (\eta)) = 0.
	\end{equation}

	Looking at the definition of the evolution vector field~\eqref{evo-def}, we obtain the following result.
	\begin{corollary}\label{corola}
		The map	$(X,f):\mathcal{M}\mapsto \mathcal{T}\mathcal{M}$ defines a Lagrangian submanifold of the exact symplectic manifold $(H\mathcal{M}\times \mathbb{R}, \omega_{\eta})$ if and only if, locally, $X = \varepsilon_H$ and $f = \mathcal{R}(H)$ for a (local) smooth function $H: \mathcal{M} \to \mathcal{R}$.
	\end{corollary}
	\textbf{Proof.}
	Since the evolution vector field satisfies
	\begin{equation}
		\mathcal{L}_{\varepsilon_H}\eta=dH-\mathcal{R}(H)\eta,\qquad \eta(\varepsilon_H)=0,
	\end{equation}
one arrives that $\mathcal{L}_X \eta + f \eta = d H$ is closed.
	Conversely, if $\eta(X) = 0$ and $\mathcal{L}_X \eta + f \eta$ is closed, in a local chart $U$, one has that  
	\begin{equation}
		\left.\mathcal{L}_X \eta + f \eta\right \vert_U = \left. \contr{X} d \eta + f \eta\right \vert_U  = d H 
	\end{equation}
	for a local smooth  function $H$. By contracting the equality with the Reeb vector field, we obtain that $f = \mathcal{R}(H)$. Thus, $X = \varepsilon_H$.
	
\textbf{The right wing of the triple.}
We start by presenting the following commutative diagram:
\begin{equation}
	\begin{tikzcd}
		{(\mathcal{T}\mathcal{M},\eta^{\mathcal{T}})} \arrow[rr, "\beta^c"]                      &                                                                      & {(\mathcal{T}^*\mathcal{M}, \eta_\mathcal{M})}            \\ \\
		{(H \mathcal{M}  \times \mathbb{R}, \omega_{\eta})} \arrow[uu, "j", hook] \arrow[rr, "\beta^0"] &                                                                      & {(T^*\mathcal{M}, \omega_\mathcal{M})}. \arrow[uu, "{(\mathrm{Id}_{T\mathcal{M}},0)}",swap, hook] 		\end{tikzcd}
	\end{equation}
	Here,  $\beta^c$ is the contact mapping in \eqref{beta-c}, the mapping $(\mathrm{Id}_{T\mathcal{M}},0)$ is the canonical inclusion of the cotangent bundle into the extended cotangent bundle as a zero section and the smooth map $\beta^0$ is given by
	\begin{equation}
	 \beta^0:H \mathcal{M}  \times \mathbb{R}\longrightarrow T^*\mathcal{M},\qquad (V,u)\mapsto -\contr{V} d \eta - u \eta.
	 \end{equation}
In coordinates, 
\begin{equation}
\beta^0(q^i,p_i,z,\dot{q}^i,\dot{p}_i, u) =(q^i,p_i,z, u p_i+\dot{p}_i ,- \dot{q}^i, -u).
	 \end{equation}

Now we are ready to plot the right wing of Tulczyjew's triple for the evolution  contact dynamics. For this, we replace the contact manifold $\mathcal{M}$ with the extended cotangent bundle $\mathcal{T}^*\mathcal{Q}$. Note that, in this case, we have
\begin{equation}\label{Right-w-evo}
	\begin{tikzcd}
		 H\mathcal{T}^*\mathcal{Q}   \times \mathbb{R} \arrow[rr, "\beta^0"] \arrow[rdd,"\hat{\tau}_{\mathcal{T}^*\mathcal{Q}}"]&                                                                      & T^*\mathcal{T}^*\mathcal{Q} \arrow[ldd,swap,"\pi_{\mathcal{T}^*\mathcal{Q}}"] \\  \\
										  & \mathcal{T}^*\mathcal{Q} \arrow[luu, bend left, "{(\varepsilon_H, \mathcal{R}(H))}"] \arrow[ruu, bend right, swap,"-d H"] &                                                             
		\end{tikzcd}
	\end{equation}
where $\hat{\tau}_{\mathcal{T}^*\mathcal{Q}}$ is the projection taking a two-tuple $(V,u)$ mapping to $\tau_{\mathcal{T}^*\mathcal{Q}}(V)$, that is,
 \begin{equation}
 \hat{\tau}_{\mathcal{T}^*\mathcal{Q}}: H\mathcal{T}^*\mathcal{Q}   \times \mathbb{R}\longrightarrow \mathcal{T}^*\mathcal{Q},\qquad (V,u)\mapsto \tau_{\mathcal{T}^*\mathcal{Q}}(V)
 \end{equation}
using that $ H\mathcal{T}^*\mathcal{Q} $ is a vector subbundle of $ T\mathcal{T}^*\mathcal{Q} $. In order to see that the triangle commutes, we compute that
\begin{equation}
	\beta^0\circ (\varepsilon_H,\mathcal{R}(H)) = - \contr{\varepsilon_H} \dd \eta_\mathcal{Q} - \mathcal{R}(H) \eta_\mathcal{Q} = -\mathcal{L}_{\varepsilon_H} \eta_\mathcal{Q}- \mathcal{R}(H) \eta_\mathcal{Q} = - d H,
\end{equation}
where we have used~\eqref{evo-def} and the Cartan's formula. Following notation presented in \eqref{sss}, we write that the quintuple 
\begin{equation}  
		( H\mathcal{T}^*\mathcal{Q}   \times \mathbb{R},\hat{\tau}_{\mathcal{T}^*\mathcal{Q}},\mathcal{T}^*\mathcal{Q},-\theta_{\eta_{\mathcal{Q}}},\beta^0)
		\end{equation}
		determines a special symplectic structure. 
Since $\beta^0$ is a symplectic diffeomorphism, we can realize one  more time that the image space $im(\varepsilon_H,\mathcal{R}(H))$ as a Lagrangian submanifold of $H\mathcal{T}^*\mathcal{Q}   \times \mathbb{R}$. 
Moreover, if $\sigma^e:I\mapsto \mathcal{T}^*\mathcal{Q}$ is a smooth curve on $\mathcal{T}^*\mathcal{Q}$ and we consider the lift $(\sigma^e)^{\mathcal{T}}:I\mapsto \mathcal{T}\mathcal{T}^*\mathcal{Q}$ of the curve $\sigma^e$ given by \eqref{sigma-T}, then one may prove that $\sigma^e$ is a solution of the contact evolution equations for $H$ if and only if $(\sigma^e)^{\mathcal{T}}(t)\in im(\varepsilon_H,\mathcal{R}(H))$, for every $t\in I$.

\textbf{The left wing part of the triple.}
Recall the mapping $\alpha^c$ given in \eqref{alpha-c}. Consider the following commutative diagram
\begin{equation}\label{comdiag-TTT}
	\begin{tikzcd}
		{(\mathcal{T} \mathcal{T}^*\mathcal{Q},\eta^{\mathcal{T}})} \arrow[rr, "\alpha^c"]   &                                                 & {(\mathcal{T}^* \mathcal{T}\mathcal{Q}, \eta_{ \mathcal{T} \mathcal{Q}})} \arrow[dd, "pr",swap]  \\\\
		{(H\mathcal{T}^*\mathcal{Q} \times \mathbb{R}, \omega_{\eta_{\mathcal{Q}}})} \arrow[uu, "j", hook] \arrow[rr,"{\alpha^0}"] &                                                 & {(T^*\mathcal{T}\mathcal{Q}, \omega_{\mathcal{T}\mathcal{Q}})}  
		\end{tikzcd}
\end{equation}
where $j$ is the inclusion mapping in \eqref{map-j},  $pr$ is the natural projection  from the extended cotangent bundle $\mathcal{T}^* \mathcal{T}\mathcal{Q}=T^* \mathcal{T}\mathcal{Q}\times \mathbb{R}$ to the first factor that is the cotangent bundle $T^* \mathcal{T}\mathcal{Q}$. Using  
\eqref{alpha-c} and \eqref{map-j}, we deduce that the local expression of 
 $\alpha^0$ is
\begin{equation}\label{alfa-0}
\alpha^0: H\mathcal{T}^*\mathcal{Q} \times \mathbb{R}\longrightarrow T^*\mathcal{T}\mathcal{Q},\qquad 
(q^i,p_i,z,\dot{q}^i,\dot{p}_{i}, u)\mapsto 
(q^i,\dot{q}^i,z,up_i+\dot{p}_{i},p_i,-u).
 \end{equation}
A direct calculation proves that  $(\alpha^0)^*\omega_{\mathcal{T}\mathcal{Q}}=\omega_{\eta_{\mathcal{Q}}}$ that is $\alpha^0$ is a symplectic diffeomorphism. Further, we can determine another potential one-form $\theta_{\eta_{\mathcal{Q}}}'$ for the symplectic two-form $\omega_{\eta_{\mathcal{Q}}}$ given in \eqref{omega-H} as follows
\begin{equation}
\theta_{\eta_{\mathcal{Q}}}':=(\alpha^0)^*\theta_{\mathcal{T}\mathcal{Q}}=up_idq^i+\dot{p}_id\dot{q}^i-udz,
 \end{equation}
where $\theta_{\mathcal{T}\mathcal{Q}}$ is the canonical one-form on the cotangent bundle $T^*\mathcal{T}\mathcal{Q}$. 
We now construct the left wing of evolution  contact Tulczyjew's triple as follows
\begin{equation}\label{Left-w-evo}
	\begin{tikzcd}
		 H\mathcal{T}^*\mathcal{Q}   \times \mathbb{R} \arrow[rr, "\alpha^0"] \arrow[rdd,"\widetilde{\mathcal{T}\pi_{\mathcal{Q}}^0}"]&                                                                      & T^*\mathcal{T} \mathcal{Q} \arrow[ldd,swap,"\pi_{\mathcal{T}\mathcal{Q}}"] \\  \\
											  & \mathcal{T}\mathcal{Q}  \arrow[ruu, bend right, swap,"dL"]   &                                                             
		\end{tikzcd}
	\end{equation}
 where $\widetilde{\mathcal{T}\pi_{\mathcal{Q}}^0}$ is the projection defined by $\widetilde{\mathcal{T}\pi_{\mathcal{Q}}^0}= \mathcal{T}\pi_{\mathcal{Q}}^0 \circ j=\pi_{\mathcal{T}\mathcal{Q}}\circ \alpha^0$. The local expression of $\widetilde{\mathcal{T}\pi_{\mathcal{Q}}^0}$ is
 \begin{equation}\label{tilde-map}
\widetilde{\mathcal{T}\pi_{\mathcal{Q}}^0}: H\mathcal{T}^*\mathcal{Q}\times \mathbb{R}\longrightarrow \mathcal{T}\mathcal{Q},\qquad (q^i,p_i,z,\dot{q}^i,\dot{p}_{i}, u)\mapsto (q^i,\dot{q}^i,z).
 \end{equation}
Note that, following the notation in \eqref{sss}, we have that the quintuple
 \begin{equation}
(H\mathcal{T}^*\mathcal{Q}\times \mathbb{R},\widetilde{\mathcal{T}\pi_{\mathcal{Q}}^0},\mathcal{T}\mathcal{Q},-\theta_{\eta_\mathcal{Q}}',\alpha^0)  
 \end{equation} 
 is a special symplectic structure.

For a given Lagrangian function $L:\mathcal{T}\mathcal{Q}\mapsto \mathbb{R}$, $\alpha^0$ pulls the Lagrangian submanifold $im(dL)$ to a Lagrangian submanifold $(\alpha^0)^{-1}(im(dL))$ of  the symplectic manifold $H\mathcal{T}^*\mathcal{Q} \times \mathbb{R}$. This  Lagrangian submanifold is the realization of the evolution Herglotz equations. In fact, 
 \begin{equation}
 (\alpha^0)^{-1}(im(dL))=\big\{(q^i,\frac{\partial L}{\partial \dot{q}^i},z,\dot{q}^i,\frac{\partial L}{\partial z}\frac{\partial L}{\partial \dot{q}^i}+\frac{\partial L}{\partial q^i},\dot{z}, -\frac{\partial L}{\partial z})  \in 
 \mathcal{T} \mathcal{T}^*\mathcal{Q}:\dot{z}=\dot{q}^i\frac{\partial L}{\partial \dot{q}^i}
  \big\}.
  \end{equation} 
So, a curve $c^e:I\mapsto \mathcal{T}\mathcal{Q}$ is a solution of the evolution Herglotz equations for $L$ if and only if $(\mathbb{F}L\circ c^e)^\mathcal{T}(t) \in  (\alpha^0)^{-1}(im(dL))$, for every $t\in I$, where  $(\mathbb{F}L\circ c^e)^\mathcal{T}:I\mapsto \mathcal{T}\mathcal{T}^*\mathcal{Q}$ is the lift of the curve $\mathbb{F}L\circ c^e:I\mapsto \mathcal{T}^*\mathcal{Q} $ given by \eqref{sigma-T}.

We merge the right and left wings of triple given in \eqref{Right-w-evo} and \eqref{Left-w-evo} respectively. So, we have the evolution  contact Tulczyjew's triple:
\begin{equation}\label{TT-evo}
	\begin{tikzcd}
	 T^*\mathcal{T} \mathcal{Q}	 \arrow[rdd,"\pi_{\mathcal{T}\mathcal{Q}}"]&                                                                      & H\mathcal{T}^*\mathcal{Q}   \times \mathbb{R}  \arrow[ldd,swap,"\widetilde{\mathcal{T}\pi_{\mathcal{Q}}^0}"] \arrow[ll, "\alpha^0",swap] \arrow[rr, "\beta^0"] \arrow[rdd,"\hat{\tau}_{\mathcal{T}^*\mathcal{Q}}"]&&                                                                     T^*\mathcal{T}^*\mathcal{Q} \arrow[ldd,swap,"\pi_{\mathcal{T}^*\mathcal{Q}}"]\\  \\
											  & \mathcal{T}\mathcal{Q} \arrow[rdd,"\tau^0_{\mathcal{Q}}",swap]   && \mathcal{T}^*\mathcal{Q} \arrow[ldd,"\pi^0_{\mathcal{Q}}"]
										\\	  \\&&\mathcal{Q}
		\end{tikzcd}
	\end{equation}
Note that this triple is consisting of two special symplectic structures.

\subsection{The Legendre Transformation}\label{Sec-Leg-trf}

At first we introduce the following Whitney product
\begin{equation}\label{WP-2}
\mathcal{T} \mathcal{Q} \times_{\mathcal{Q}\times \mathbb{R}} \mathcal{T}^*\mathcal{Q}=\{(u,\zeta,z)\in T\mathcal{Q} \times T^*\mathcal{Q}\times \mathbb{R}: \tau_{\mathcal{Q}}(u)= \pi_{\mathcal{Q}}(\zeta)\}
\end{equation} 
over the extended manifold $\mathcal{Q}\times \mathbb{R}$. This is an $3n+1$ dimensional manifold which we call extended Pontryagin bundle. We introduce the local coordinates $(q^i,\dot{q}^i,p_i,z)$ on $\mathcal{T} \mathcal{Q} \times_{\mathcal{Q}\times \mathbb{R}} \mathcal{T}^*\mathcal{Q}$. 
We now 	generate the Legendrian submanifold $\mathcal{N}_L$ on $\mathcal{T}\mathcal{T}^*\mathcal{Q}$ given in \eqref{N_L} referring to  the right wing of the contact Tulczyjew's triple (\ref{con-TT}). This is the Legendre transformation in the understanding of Tulczyjew. To have that, define the following Morse family
\begin{equation} \label{energy}
E:\mathcal{T} \mathcal{Q} \times_{\mathcal{Q}\times \mathbb{R}} \mathcal{T}^*\mathcal{Q}\longrightarrow \mathbb{R},\qquad (q^i,\dot{q}^i,p_i,z)\mapsto p_i \dot{q}^i-L(q,\dot{q},z).
\end{equation}
To represent the minus of the Morse family $-E$, we are drawing the right side of the Tulczyjew's triple (\ref{con-TT}) by equipping it with the Whitney product \eqref{WP-2} as follows
\begin{equation}\label{Ham-Morse-Gen-2}
\xymatrix{\mathcal{T} \mathcal{T}^* \mathcal{Q} \ar[rr]^{ \beta^c} \ar[ddr]_{\tau^0_{\mathcal{T}^* \mathcal{Q}}}
&&\mathcal{T}^* \mathcal{T}^*\mathcal{Q}\ar[ddl]^{\pi^0_{\mathcal{T}^*\mathcal{Q}}} &\mathcal{T} \mathcal{Q} \times_{\mathcal{Q}\times \mathbb{R}} \mathcal{T}^*\mathcal{Q} \ar[dd]^{{\rm pr}_2}\ar[rr]^{\qquad -E} &&  \mathbb{R}
\\ \\ &\mathcal{T}^*\mathcal{Q}  \ar@{=}[rr]&&\mathcal{T}^* \mathcal{Q}}.
\end{equation}
This geometry fits well in the general picture in \eqref{Morse-Gen--}. Recalling \eqref{MFGen-C}, a direct computation determines the Legendrian submanifold of $\mathcal{T}^*\mathcal{T}^*\mathcal{Q}$ generated by  minus of the Morse family $-E$ 
as 
\begin{equation} \label{LagSubEL-T*T*--}
\begin{split}
\mathcal{N}_{-E}&=
\big\{(q^i,p_i,z, -\frac{\partial E}{\partial q^i},-\frac{\partial E}{\partial p_i},-\frac{\partial E}{\partial z},-E )\in \mathcal{T}^*\mathcal{T}^*\mathcal{Q}: \frac{\partial E}{\partial \dot{q}^i}=0\big\}
\\&=\big\{(q^i,p_i,z, \frac{\partial L}{\partial q^i},-\dot{q}^i ,\frac{\partial L}{\partial z},-p_i \dot{q}^i+L)\in \mathcal{T}^*\mathcal{T}^*\mathcal{Q}: p_i-\frac{\partial L}{\partial \dot{q}^i}=0\big\}.
\end{split}
\end{equation}
Using the inverse of the contact diffeomorphism $\beta^c$, we map the Legendrian submanifold $\mathcal{N}_{-E}$ to a Legendrian submanifold of  the contact manifold $\mathcal{T}\mathcal{T}^*\mathcal{Q}$ as 
\begin{equation} \label{LagSubEL-22}
(\beta^c)^{-1}(\mathcal{N}_{-E})=\Big\{\big(q^i,p_i,z, \dot{q}^i, p_i\frac{\partial L}{\partial z}+ \frac{\partial L}{\partial q^i}, L, -\frac{\partial L}{\partial z}\big) \in \mathcal{T}\mathcal{T}^*\mathcal{Q}: p_i-\frac{\partial L}{\partial \dot{q}^i}=0\Big\}
\end{equation}
See that, this Legendrian submanifold is exactly the Legendrian submanifold $\mathcal{N}_L$ in \eqref{N_L} realizing the Herglotz equations \eqref{Herglotz} and \eqref{herglotzprinciple}. This completes the Legendre transformation of the Herglotz equations to the contact Hamiltonian formalism.

If, further, the Lagrangian function $L=L(q,\dot{q},z)$ is non-degenerate  then from the equation
\begin{equation} \label{FL}
\frac{\partial E}{\partial \dot{q}^i}(q,\dot{q},p,z)=p_i-\frac{\partial L}{\partial\dot{q}^i}(q,\dot{q},z)=0
\end{equation} 
one can explicitly determine the velocity $\dot{q}^i$ in terms of $(q^i,p_i,z)$. This gives
\begin{equation} \label{FL-}
\mathbb{F}L:\mathcal{T}\mathcal{Q}\longrightarrow \mathcal{T}^*\mathcal{Q},\qquad  (q^i ,\dot{q}^j,z) \longrightarrow \big( q^i ,\frac{\partial L}{\partial \dot{q}^j}( q ,\dot{q},z),z\big)
\end{equation}
as a local diffeomorphism relating $\mathcal{T}\mathcal{Q}$ and $\mathcal{T}^*\mathcal{Q}$. In this case, the Morse family $E$ can be reduced to a well-defined contact Hamiltonian function
\begin{equation}
H(q ,p ,z)=p_i \dot{q}^i( q,p,z)  -L\big( q ,\dot{q}( q ,p,z),z\big)
\end{equation}
on $\mathcal{T}^*\mathcal{Q}$.

It is possible to perform the inverse Legendre transformation of the contact Hamiltonian dynamics as well. This time, one needs to generate the Legendrian submanifold $\mathcal
{N}_{-H}$ in \eqref{N_H} realizing the contact Hamilton's equations \eqref{conHamee} referring to the left wing of the contact triple \eqref{con-TT}. 

\textbf{The Legendre Transformation for Evolution  Dynamics.} Recall the Tulczyjew's triple \eqref{TT-evo} exhibited for the case of evolution  contact dynamics. 
We consider the Lagrangian submanifold $(\alpha^0)^*(im(dL))$ of $H\mathcal{T}^*\mathcal{Q}\times \mathbb{R}$ generated by a Lagrangian $L=L(q,\dot{q},z)$ on  $\mathcal{T}\mathcal{Q}$ referring to the left wing \eqref{Left-w-evo} of the triple. Once more consider the total space given in \eqref{WP-2}, and the energy function $E$ given in \eqref{energy}. In this evolution  case, we plot the following diagram merging the right wing \eqref{Right-w-evo} of the evolution  Tulczyjew's triple and the Morse family determined by $-E$ that is 
\begin{equation}\label{Ham-Morse-Gen-3}
\xymatrix{H\mathcal{T}^*\mathcal{Q}\times \mathbb{R} \ar[rr]^{ \beta^0} \ar[ddr]_{ \hat{\tau}_{\mathcal{T}^*\mathcal{Q}}}
&&T^* \mathcal{T}^*\mathcal{Q}\ar[ddl]^{\pi_{\mathcal{T}^*\mathcal{Q}}} &\mathcal{T} \mathcal{Q} \times_{\mathcal{Q}\times \mathbb{R}} \mathcal{T}^*\mathcal{Q} \ar[dd]^{{\rm pr}_2}\ar[rr]^{\qquad -E} &&  \mathbb{R}
\\ \\ &\mathcal{T}^*\mathcal{Q}  \ar@{=}[rr]&&\mathcal{T}^* \mathcal{Q}}.
\end{equation}
From \eqref{MFGen}, we deduce that the Lagrangian submanifold $\mathcal{S}_{-E}$ of the cotangent bundle $T^*\mathcal{T}^*\mathcal{Q}$ generated by $-E$ is computed to be
\begin{equation}
\begin{split}
\mathcal{S}_{-E}&=
\big\{(q^i,p_i,z, -\frac{\partial E}{\partial q^i},-\frac{\partial E}{\partial p_i},-\frac{\partial E}{\partial z})\in T^*\mathcal{T}^*\mathcal{Q}: \frac{\partial E}{\partial \dot{q}^i}=0\big\}
\\&=\big\{(q^i,p_i,z, \frac{\partial L}{\partial q^i},-\dot{q}^i ,\frac{\partial L}{\partial z})\in \mathcal{T}^*\mathcal{T}^*\mathcal{Q}: p_i-\frac{\partial L}{\partial \dot{q}^i}=0\big\}.
\end{split}
\end{equation}
Using the inverse of the symplectic diffeomorphism $\beta^0$, we transfer the Lagrangian submanifold $\mathcal{S}_{-E}$ to a Lagrangian submanifold of $H\mathcal{T}^*\mathcal{Q}\times \mathbb{R}$ as follows
\begin{equation*} 
(\beta^0)^{-1}(\mathcal{S}_{-E})=\Big\{\big(q^i,p_i,z, \dot{q}^i, p_i\frac{\partial L}{\partial z}+ \frac{\partial L}{\partial q^i}, \dot{z}, -\frac{\partial L}{\partial z}\big) \in H\mathcal{T}^*\mathcal{Q}\times \mathbb{R}: p_i-\frac{\partial L}{\partial \dot{q}^i}=0,\dot{z}-\dot{q}^i\frac{\partial L}{\partial \dot{q}^i}=0 \Big\}.
\end{equation*}
This is exactly the Lagrangian submanifold $(\alpha^0)^{-1}(im(dL))$ realizing the evolution  Herglotz equations. 
In a similar way, one may obtain the inverse Legendre transformation of the contact evolution dynamics for a Hamiltonian function $H:\mathcal{T}^*\mathcal{Q}\mapsto \mathbb{R}$. 

\section{Example: The Ideal Gas}

\subsection{A Quantomorphism on the Euclidean Space}

Thermodynamics has been studied extensively in the framework of contact geometry. For some recent work directly related with the present discussions, we cite  \cite{Br17,bravetti2019contact,ghosh2019contact,
mrugala1991contact}. In this section, we shall be applying the theoretical results obtained in the previous sections to some thermodynamical models. 

We start this subsection by providing the following theorem realizing a 
strict contact diffeomorphism (quantomorphism) on the extended cotangent bundle $\mathcal{T}^*\mathbb{R}^{m}$ of the Euclidean space, see also \cite{bravetti2019contact}. The proof  follows by a direct calculation. 
  \begin{theorem}\label{theo-naive}
Consider a disjoint partition $I \cup J$ of the set of indices $\{1,\dots, m\}$ so that the coordinates on $\mathbb{R}^{m}$ is given as $(x^a,x^\rho)$, where $a\in I$ and $\rho\in J$. Then the following  mapping
\begin{equation}\label{phi}
\phi:\mathcal{T}^*\mathbb{R}^m\longrightarrow \mathcal{T}^*\mathbb{R}^m,\qquad (x^a,x^\rho,y_a,y_\rho,u)\mapsto (x^a,y_\rho,y_a,-x^\rho,u-x^\rho y_\rho) 
\end{equation}
preserves the canonical contact one-form $\eta_{\mathbb{R}^{m}}=du-y_adx^a-y_\rho d x^\rho$. Here, $(x^a,x^\rho,y_a,y_\rho,u)$ are the Darboux' coordinates on the extended cotangent bundle $\mathcal{T}^*\mathbb{R}^m$. 
\end{theorem}

In Subsection \ref{sec-sub}, we have stated that the image of the first prolongation of a smooth function on the base manifold is a Legendrian submanifold of the extended cotangent bundle. Accordingly, consider a smooth function $U=U(x^a,x^\rho)$ on $\mathbb{R}^m$ so that its first prolongation $\mathcal{T}^*U$ to the extended cotangent bundle turns out to be a Legendrian submanifold of $\mathcal{T}^*\mathbb{R}^{m}$ as given in \eqref{j1F}. Under the quantomorphism  $\phi$ in \eqref{phi}, we have a Legendrian submanifold on the image space as
\begin{equation}\label{phi-2}
\phi(x^a,x^\rho,\frac{\partial U}{\partial x^a},\frac{\partial U}{\partial x^\rho},U)=(x^a,\frac{\partial U}{\partial x^\rho}, 
\frac{\partial U}{\partial x^a},-x^\rho,U-x^\rho\frac{\partial U}{\partial x^\rho}).
\end{equation}
This alternative realization of the Legendrian submanifold is important for geometric characterization of reversible thermodynamics.   

\begin{remark}
If $\mathcal{Q}=\mathbb{R}^n$ then the extended cotangent and the extended tangent bundles turn out to be isomorphic that is $\mathcal{T}^*\mathcal{Q}\simeq \mathcal{T}\mathcal{Q}\simeq \mathbb{R}^{2n+1}$. In this particular instance, by assuming $m=2n+1$, the mapping $\psi^c$ in \eqref{psi-c} can be regarded as a particular case of $\phi$ in \eqref{phi} if the canonical coordinates $(q^i,p_i,z)$ are decomposed as $(x^a)=(q^i,z)$ where $a=1,\dots ,n+1$, and $(x^\rho)=(p_i)$ where $\rho=1,\dots ,n$. 
\end{remark}

\subsection{Equilibrium Thermodynamics} 

Obeying the geometry exhibited in the previous section, we take $m=3$ with coordinates $(S,V,N)\in \mathbb{R}^3$. Here, $S$ stands for the entropy, $V$ is the volume, and $N$ is the mole number of classical ideal gas. The conjugate variables $(T,-P,\mu)\in (\mathbb{R}^3)^*$ are the  temperature, the pressure, and the chemical potential, respectively. By employing the internal energy $U$ as the fiber coordinate, we complete the following realization of the extended cotangent bundle $(S,V,N,T,-P,\mu,U)\in \mathcal{T}^*\mathbb{R}^{3}$. 
Consider the contact one-form 
\begin{equation}
 \eta_{\mathbb{R}^3}=dU-TdS+PdV-\mu dN.
\end{equation}

As a particular instance, we choose the internal energy\begin{equation}\label{int-en}
U(S,V,N)=U_0V^{-1/c}N^{(c+1)/c}\exp(\frac{S}{cNR}),
\end{equation} 
as a function depending on the base coordinates $(S,V,N)\in \mathbb{R}^3$. Here, $U_0$ is a positive constant, $c$ is the heat capacity and $R$ is the universal gas constant. The first prolongation $\mathcal{T}^*U$ is a Legendrian submanifold $\mathcal{N}$ of the contact manifold  $(\mathcal{T}^*\mathbb{R}^{3},\eta_{\mathbb{R}^3})$. By considering that the temperature $T=\partial U/ \partial S$ and the pressure $P=-\partial U/ \partial V$, we have the following set of  equations 
\begin{equation}\label{gas}
cV^{1/c}RT=U_0N^{1/c}\exp(\frac{S}{cNR}), \qquad PV=NRT, \qquad \mu=(c+1)RT-TS/N
\end{equation}
those realizing $\mathcal{N}$.

\textbf{The Legendre Transformations.} In the light of Theorem \ref{theo-naive}, and the transformation \eqref{phi-2}, we now present the Legendre transformation between the internal energy, the enthalpy, the Helmholtz function, and the Gibbs function. For a similar discussion but in the framework of symplectic geometry see \cite{Tu77}. We start with the Legendrian submanifold determined by the internal energy $U$ in \eqref{int-en}. 
\newline
\textbf{(1)} We decompose the base variables as $(S,N)$ and $V$ and apply Theorem \ref{theo-naive} to the volume variable. This is resulting with a quantomorphism  computed to be
 \begin{equation}  
 \phi^1:\mathcal{T}^*\mathbb{R}^3\longrightarrow \mathcal{T}^*\mathbb{R}^3,\qquad (S,V,N,T,-P,\mu,U)\mapsto (S,-P,N,T,-V,\mu,U+PV).
 \end{equation}
Note that, on the image space, the fiber component is the enthalpy function $B=U+PV$. If we solve the pressure from the equation $P=-\partial U/ \partial V$, the enthalpy function can be written as a function of the new base variables $(S,P,N)$ that is 
 \begin{equation}  \label{entalphy}
 B(S,P,N)=\bar{c}~U_0^{(c+1)/c}P^{1/(1+c)}N\exp(\frac{S}{(c+1)NR}),
  \end{equation}
  where $\bar{c}$ is a constant defined to be $c^{1/(1+c)}+c^{-c/(1+c)}$.
So that, the enthalpy function is another generator of the same Legendrian submanifold. Indeed, the first prolongation of $\mathcal{T}^*B$ is given the system of equations in \eqref{gas} so that $ \phi^1 \circ \mathcal{T}^*B=\mathcal{N}$. \newline
 \textbf{(2)} We start once more with the internal energy but this time we perform the transformation to the entropy variable  $S$. For this case, we have the quantomorphism 
  \begin{equation} \label{hoo-2} 
 \phi^2:\mathcal{T}^*\mathbb{R}^3\longrightarrow \mathcal{T}^*\mathbb{R}^3,\qquad (S,V,N,T,-P,\mu,U)\mapsto (T,V,N,-S,-P,\mu,U-ST).
 \end{equation}
In this case, the fiber term $F=U-ST$ is the Helmholtz function. Using the identity $T=\partial U /\partial S$, we write the Helmholtz function as a function of the base components $(T,V,N)$ of the image space, that is
 \begin{equation}  \label{Helmholtz}
 F(T,V,N)=cNRT\big(1+\frac{1}{c}\log N - \frac{1}{c}\log V +\log (\frac{U_0}{cRT})\big ).
  \end{equation}
So that, $F$ is another generator of the same Legendrian submanifold determined by the equations \eqref{gas}, that is $\phi^2 \circ \mathcal{T}^*F=\mathcal{N}$. \newline
\textbf{(3)} This time, we consider the Helmholtz function $F$ in \eqref{Helmholtz} and apply the transformation given in Theorem \ref{theo-naive} to the volume variable, that is
  \begin{equation} \label{hoo-3} 
  \phi^3: \mathcal{T}^*\mathbb{R}^3\longrightarrow \mathcal{T}^*\mathbb{R}^3,\qquad  (T,V,N,-S,-P,\mu,F)\mapsto (T,-P,N,-S,-V,\mu,F+PV).
   \end{equation}
The fiber term $G=F+PV$ is the Gibbs function. By taking $P=-\partial F / \partial V$, we write the Gibbs function as a function of the base variables, that is
  \begin{equation} \label{Gibbs}
  G(T,P,N)=NRT\big(1+c+\log N - \log \frac{NRT}{P} +c\log (\frac{U_0}{cRT})\big ).
     \end{equation}
We have that $\phi^3 \circ \mathcal{T}^*G=\mathcal{N}$ is the Legendrian submanifold determined by the equations \eqref{gas}. Evidently, by iteratively applying $\phi^2$ in \eqref{hoo-2} and $\phi^3$ in \eqref{hoo-3}  one can define a quantomorphism from the internal energy setting to the  Gibbs function setting. In this case the quantomorphism is determined as 
  \begin{equation} \label{quantato}
  \phi^3\circ   \phi^2 :\mathcal{T}^*\mathbb{R}^3\longrightarrow \mathcal{T}^*\mathbb{R}^3, \qquad (S,V,N,T,-P,U)\mapsto (T,-P,N,-S,-V,\mu,U-TS+PV).
       \end{equation}
This gives  a direct passage from the internal energy to the Gibbs function. 
\newline
\textbf{(4)} To complete the Legendre transformations, we consider the Gibbs function in \eqref{Gibbs} and define the quantomorphism to the chemical potential variable $\mu$, that is
  \begin{equation} \label{phi-4}
  \phi^4: \mathcal{T}^*\mathbb{R}^3\longrightarrow \mathcal{T}^*\mathbb{R}^3,\qquad  (T,-P,N,-S,-V,\mu,G)\mapsto (T,-P,\mu,-S,-V,-N,G-\mu N).
   \end{equation}
   The fiber variable $W=G-\mu N$ determines a new generator of the Legendrian submanifold $\mathcal{N}$ depending on the base variables of the image space. In this case we have that 
     \begin{equation} \label{Wibbs}
  W(T,P,\mu)=NRT\big(1+c+ \log \frac{ST}{(c+1)RT-\mu} - \log \frac{NRT}{P} +c\log (\frac{U_0}{cRT})\big )-\frac{\mu ST}{(c+1)RT-\mu}.
     \end{equation}
Composing the quantomorphism $\phi^3\circ   \phi^2 $ in \eqref{quantato} with $\phi^4$ in \eqref{phi-4} one can derive a Legendre transformation from the internal energy $U$ setting to $W$ as
  \begin{equation} \label{quantato-2}
 \phi :\mathcal{T}^*\mathbb{R}^3\longrightarrow \mathcal{T}^*\mathbb{R}^3, \quad (S,V,N,T,-P,\mu,U)\mapsto (T,-P,\mu,-S,-V,-N,U-TS+PV-\mu N).
       \end{equation}
This is the full Legendre transformation of the ideal gas. If we insists that the dynamics of the gas stays in the Legendrian submanifold $\mathcal{N}$ determined through the equations \eqref{gas}, we can consider the fiber variable $U-TS+PV-\mu N$ for an arbitrary internal energy $U$ as a Hamiltonian function of the dynamics for the ideal case.  
Let us discuss this motion in the following subsection. 

\subsection{Hamiltonian Flow and Its Legendrian Realization}
 We consider the  quantomorphism in \eqref{quantato-2} but instead of a specific internal energy, we let $U$ be an independent variable. Consider the Hamiltonian function
\begin{equation}\label{Ham-gas}
H(S,V,N,T,-P,\mu,U)=TS-NRT+\mu N-U
\end{equation}  
on the extended cotangent bundle $\mathcal{T}^*\mathbb{R}^{3}$. Note that, we determine the Hamiltonian function \eqref{Ham-gas}  by substituting the ideal gas equation $PV=NRT$ (that is the second equation in \eqref{gas}) into the minus of the generator function $U-TS+PV-\mu N$ of the quantomorphism. The minus sign is to fit physical intuition, for example, to be sure that the entropy $S$ is increasing along the motion. In fact, while computing the Hamiltonian dynamics, the minus sign will be compensated with the minus sign appearing in the formulation given in \eqref{contact-Ham-Leg}. According to \eqref{Ham-Leg-Sub}, the image space of the first prolongation of the Hamiltonian function \eqref{Ham-gas}
\begin{equation}\label{ex-Leg-1}
\im(-\mathcal{T}^*H)=\big\{ (S,V,N,T,-P,\mu,U,-T,0,RT-\mu,-S+NR,0,-N,1,-H)\in \mathcal{T}^*\mathcal{T}^*\mathbb{R}^{3}
\big\}
\end{equation} 
is a Legendrian submanifold of the iterated extended cotangent bundle $\mathcal{T}^*\mathcal{T}^*\mathbb{R}^{3}$. On the other hand, the contact Hamiltonian vector field associated with the Hamiltonian function \eqref{Ham-gas}
 is
\begin{equation} 
X_H^c =(S-NR) \frac{\partial }{\partial S} +N \frac{\partial }{\partial N}
+P \frac{\partial }{\partial P}+ RT \frac{\partial }{\partial \mu}+U \frac{\partial }{\partial U}.
\end{equation} 
Notice that, the coefficients of $\partial/\partial T$ and $\partial/\partial V$ are zero. This gives that the Hamiltonian dynamics is isothermal and isochoric. 

We consider the induced coordinates on the extended tangent bundle $\mathcal{T}\mathcal{T}^*\mathbb{R}^{3}$ as follows. For the base manifold we use $(S,V,N,T,-P,\mu,U)$, for the fibers of the tangent bundle $T\mathcal{T}^*\mathbb{R}^{3}$ we refer $(\dot{S},\dot{V},\dot{N},\dot{T},-\dot{P},\dot{\mu},\dot{U})$
and for the extension $\mathbb{R}$ we use $u$ as the standard coordinate on $\mathbb{R}$.  Accordingly, for the present case, the lifted contact one-form \eqref{eta-T} is computed to be
\begin{equation}
\eta^{\mathcal{T}}=d\dot{U}+udU-(\dot{T}+uT)dS+(\dot{P}+uP)dV-(\dot{\mu}+u\mu)dN-Td\dot{S}+Pd\dot{V}-\mu d \dot{N}.
\end{equation} 
The Reeb field on $\mathcal{T}^*\mathbb{R}^{3}$ is $\partial/\partial U$, and the directional derivative of the Hamiltonian function \eqref{Ham-gas} is $\mathcal{R}(H)=-1$. The image space of the coupling of the Hamiltonian vector field and $\mathcal{R}(H)$ determines a Legendrian submanifold, of the extended tangent bundle $\mathcal{T}\mathcal{T}^*\mathbb{R}^{3}$, given by
\begin{equation}\label{ex-Leg-2}
\mathcal{N}_{-H}=\im(X_H^c,\mathcal{R}(H))=\big\{ (S,V,N,T,-P,\mu,U;
S-NR,0,N,0,-P,RT,U;-1)\in \mathcal{T}\mathcal{T}^*\mathbb{R}^{3}
\big\}.
\end{equation} 
Evidently, the Legendrian submanifolds in \eqref{ex-Leg-1} and \eqref{ex-Leg-2} are related with the contactomorphism $\beta^c:\mathcal{T}\mathcal{T}^*\mathbb{R}^{3}\mapsto \mathcal{T}^*\mathcal{T}^*\mathbb{R}^{3}$ by satisfying the relation $
\beta^c\circ \mathcal{N}_{-H}=\im(-\mathcal{T}^*H)$
as we have proved in \eqref{contact-Ham-Leg}. 

Now we wish to generate the Legendrian submanifold in 
\eqref{ex-Leg-2} referring to the left wing of the contact Tulczyjew's triple \eqref{con-TT}, that is to generate it via a Lagrangian function (probably as a Morse family) defined on the extended tangent bundle $\mathcal{T}\mathbb{R}^{3}$. In the light of Subsection \ref{Sec-Leg-trf}, we now apply the inverse  Legendre transformation. The first step is to apply $\alpha^c$ in \eqref{alpha-c} to the Legendrian submanifold $\mathcal{N}_{-H}$ in \eqref{ex-Leg-2}. The image space is 
\begin{equation}\label{leg-ex-lag}
\alpha^c(\mathcal{N}_{-H})=\big\{ (S,V,N,S-NR,0,N,U;
-T,0,RT-\mu,T,-P,\mu,1;U)\in \mathcal{T}^*\mathcal{T}\mathbb{R}^{3}
\big\}
\end{equation} 
which is a Legendrian submanifold of $\mathcal{T}^*\mathcal{T}\mathbb{R}^{3}$.

We assume the coordinates $(S,V,N,\dot{S},\dot{V},\dot{N},U)$ on the extended tangent bundle  $\mathcal{T} \mathbb{R}^{3} $ and define the Whitney sum of the extended tangent and the extended cotangent bundles $\mathcal{T}\mathbb{R}^{3}\times _{\mathbb{R}^{3}\times \mathbb{R}} \mathcal{T}^* \mathbb{R}^{3}$  with coordinates $(S,V,N,\dot{S},\dot{V},\dot{N},T,-P,\mu,U)$. 
Note that, in the Whitney sum, we fix the base coordinates  $(S,V,N)$ in $\mathbb{R}^{3}$ and the extension $U$ in $\mathbb{R}$. The subscript $\mathbb{R}^{3}\times \mathbb{R}$ in the notation of Whitney sum manifests these choices. A calculation gives that the Legendrian submanifold $\alpha^c(\mathcal{N}_{-H})$ in \eqref{leg-ex-lag} is generated by a Morse family on $\mathcal{T}\mathbb{R}^{3}\times _{\mathbb{R}^{3}\times \mathbb{R}} \mathcal{T}^* \mathbb{R}^{3}$. In other words, the Lagrangian function 
 \begin{equation}\label{Lag-ex-fun}
L(S,V,N,\dot{S},\dot{V},\dot{N},U;T,-P,\mu)=T(\dot{S}-S+NR)+\mu(\dot{N}-N)+P\dot{V}+U.
\end{equation}  
 can be understood as defined on the extended tangent bundle $\mathcal{T}\mathbb{R}^{3}$ but depending on the auxiliary variables $(T,-P,\mu)$. 
 By merging the Lagrangian function 
  with the left wing of the contact Tulzyjew's triple \eqref{con-TT}, we get the following diagram
\begin{equation}\label{Ham-Morse-Gen-2-}
\xymatrix{\mathbb{R}&&\ar[ll]_{L \qquad }\mathcal{T}\mathbb{R}^{3}\times _{\mathbb{R}^{3}\times \mathbb{R}} \mathcal{T}^* \mathbb{R}^{3} \ar[dd] && \mathcal{T}^* \mathcal{T} \mathbb{R}^{3} \ar[ddr]_{\pi^0_{\mathcal{T} \mathbb{R}^{3}}} 
\ar@(ul,ur)^{
\alpha^c(\mathcal{N}_{-H})}
&&\mathcal{T} \mathcal{T}^* \mathbb{R}^{3}\ar[ll]_{\alpha^c}
\ar[ddl]^{\mathcal{T}\pi^0_{\mathbb{R}^{3}}} \ar@(ul,ur)^{ \mathcal{N}_{-H}} 
\\ \\&&
\mathcal{T} \mathbb{R}^{3} \ar@{=}[rrr]&&&\mathcal{T}\mathbb{R}^{3} }
\end{equation}
According to the local realization of the generating family in \eqref{MFGen-C}, it is a direct calculation to show that the Legendrian submanifold $\alpha^c(\mathcal{N}_{-H})$ is generated by the Morse family $L$ in \eqref{Lag-ex-fun}. 

We remark that in this formalism of thermodynamics Hamiltonians are usually singular, as in the case above, so there is not Lagrangian formulation in the classical sense. Hence, we think that Tulcyzjew triples might be a useful tool in this situation.

\subsection{Evolutionary Flow and Its Lagrangian Realization}

We once more consider the Hamiltonian function $H$ exhibited in \eqref{Ham-gas} and defined  on the extended cotangent bundle $\mathcal{T}^*\mathbb{R}^{3}$. The evolutionary vector field $\varepsilon_H$ is defined in \eqref{evo-dyn}. A direct calculation determines the evolutionary vector field for the Hamiltonian function \eqref{Ham-gas}  as 
\begin{equation} 
\varepsilon_H =(S-NR) \frac{\partial }{\partial S} +N \frac{\partial }{\partial N}
+P \frac{\partial }{\partial P}+ RT \frac{\partial }{\partial \mu}+(TS-NRT+\mu N) \frac{\partial }{\partial U}
\end{equation} 
on $\mathcal{T}^*\mathbb{R}^{3}$. By referring to Corollary \ref{corola}, we establish that the image space 
\begin{equation}\label{ex-Leg-2-2}
\im(\varepsilon_H,\mathcal{R}(H))=\big\{ (S,V,N,T,-P,\mu,U;
S-NR,0,N,0,-P,RT,TS-NRT+\mu N;-1)\in \mathcal{T}\mathcal{T}^*\mathbb{R}^{3}
\big\}
\end{equation} 
turns out to be a Lagrangian submanifold of $H\mathcal{T}^*\mathbb{R}^{3}\times \mathbb{R}$ defined in \eqref{HMxR}. Via $\alpha^0$ in \eqref{alfa-0}, we map the Lagrangian submanifold \eqref{ex-Leg-2-2} to a Lagrangian submanifold 
\begin{equation}\label{ex-Leg-2-2--}
\alpha^0\big(\im(\varepsilon_H,\mathcal{R}(H))\big)=\big\{ (S,V,N,S-NR,0,N,U;
-T,0,RT-\mu,T,-P,\mu,1)\in T^*\mathcal{T}\mathbb{R}^{3}
\big\}
\end{equation} 
in the cotangent bundle $T^*\mathcal{T}\mathbb{R}^{3}$.  
According to \eqref{MFGen}, it is immediate to see that the Lagrangian function $L$ given in \eqref{Lag-ex-fun} defined on the Whitney sum of the extended tangent and the extended cotangent bundles $\mathcal{T}\mathbb{R}^{3}\times _{\mathbb{R}^{3}\times \mathbb{R}} \mathcal{T}^* \mathbb{R}^{3}$ generates the Lagrangian submanifold in \eqref{ex-Leg-2-2--}. In order to visualize the Lagrangian submanifold, we draw the following diagram by merging the Morse family \eqref{Lag-ex-fun} and the left wing of the evolution contact Tulczyjew's triple \eqref{TT-evo} 
\begin{equation}\label{Ham-Morse-Gen-2-3-}
\xymatrix{\mathbb{R}&&\ar[ll]_{L \qquad }\mathcal{T}\mathbb{R}^{3}\times _{\mathbb{R}^{3}\times \mathbb{R}} \mathcal{T}^* \mathbb{R}^{3} \ar[dd] && T^* \mathcal{T} \mathbb{R}^{3} \ar[ddr]_{\pi_{\mathcal{T} \mathbb{R}^{3}}} 
\ar@(ul,ur)^{\alpha^0\big(\im(\varepsilon_H,\mathcal{R}(H))\big)}
&&H\mathcal{T}^* \mathbb{R}^{3}\times \mathbb{R} \ar[ll]_{\alpha^0}
\ar[ddl]^{\widetilde{\mathcal{T}\pi_{\mathbb{R}^{3}}^0}} \ar@(ul,ur)^{ \im(\varepsilon_H,\mathcal{R}(H))} 
\\ \\&&
\mathcal{T} \mathbb{R}^{3} \ar@{=}[rrr]&&&\mathcal{T}\mathbb{R}^{3} }
\end{equation}
where $\widetilde{\mathcal{T}\pi_{\mathbb{R}^{3}}^0}$ is the mapping given in \eqref{tilde-map}.
    \section{Conclusions}
    In this paper, we have used the tangent contact structure, on the extended tangent bundle $\mathcal{T}\mathcal{T}^*\mathcal{Q}$, which was introduced in \cite{IbLeMaMa97}. Referring to this, and by introducing the notion of special contact structure, we have constructed a Tulcyzjew's triple for contact manifolds, see Diagram \ref{con-TT}. This permits us to describe both the contact Lagrangian and the contact Hamiltonian dynamics as Legendrian submanifolds of $\mathcal{T}\mathcal{T}^*\mathcal{Q}$. In this formulation, the Legendre transformation is defined as a passage between two generators of the same Legendrian submanifold. Note that, this approach is free from the Hessian condition. That means, it is applicable for degenerate theories as well. We, further, present Tulcyzjew's triple for evolutionary dynamics, see Diagram \ref{TT-evo}. Instead of contact structures, the evolution triple \ref{TT-evo} is consisting of special symplectic structures. In this  construction, the contact manifold $\mathcal{T}\mathcal{T}^*\mathcal{Q}$ is substituted by the extended horizontal bundle $H\mathcal{T}^*\mathcal{Q}   \times \mathbb{R}$ which is symplectic.  We have concluded the paper by applications of the theoretical results to  geometrical foundations of some thermodynamical models.

    Here are some further questions we wish to pursue:
  \begin{itemize}
  \item In Subsection \ref{Sec-Tcont}, we have established that the image space of a contact Hamiltonian vector field is a Legendrian submanifold of the tangent contact manifold. Evidently, not all Legendrian submanifolds determine explicit dynamical equations. This observation motivates us to define the notion of an implicit Hamiltonian Contact Dynamics as a non-horizontal Legendrian submanifold of the tangent contact manifold. We refer \cite{MeMaTu95}  for a similar discussion done for the case of symplectic dynamics and integrability of the non-horizontal Lagrangian submanifolds. We find interesting to elaborate integrability of implicit Hamiltonian Contact Dynamics.  
  
  \item Following, the first question raised in this section, we plan to write  a Hamilton-Jacobi theory for implicit Hamiltonian Contact Dynamics.   
  Hamilton-Jacobi theory for (explicit) Hamiltonian Contact Dynamics is recently examined in \cite{GrPa20,de2021hamilton}. Hamilton-Jacobi theory for implicit symplectic dynamics is discussed in \cite{EsLeSa18,EsLeSa20}. 
  
  \item In the literature, Tulczyjew's triple for higher order classical dynamical systems is already available \cite{deLeLa89,esen2018geometry}. Higher order contact dynamics is studied in \cite{LeGaLaMuRo2021}. As a future work, we plan to extend the geometry presented in the present paper  to higher order contact framework.  
\end{itemize}   
    
    \section{Acknowledgment}

    	M. de Le\'on and M. Lainz acknowledge the partial finantial support from MICINN Grant PID2019-106715GB-C21 and the ICMAT Severo Ochoa project CEX2019-000904-S.  
M. Lainz wishes to thank MICINN and ICMAT for a FPI-Severo Ochoa predoctoral contract PRE2018-083203. J.C. Marrero acknowledges the partial support from European Union (Feder) grant PGC2018-098265-B-C32.

\bibliographystyle{amsplain}
\bibliography{references}
\end{document}